\documentclass[10pt]{article}
\usepackage{amsfonts,color}
 \usepackage{epsfig}
\usepackage{lscape}
\usepackage{amssymb}
 \usepackage{footnote}
\usepackage[utf8]{inputenc}
\usepackage{mathtools}
\usepackage[english]{babel}
\usepackage{hyperref}
\usepackage{cite}
\usepackage[nameinlink,capitalize]{cleveref}

\crefname{Problem}{Problem.}{Problem.}

\usepackage{graphicx}
\usepackage{bm}
\usepackage{amsfonts}
\usepackage{wrapfig}
\usepackage{subcaption}
\usepackage{amsmath}
\usepackage{float}
\usepackage{nicefrac}

\usepackage{tikz}
\usepackage{array}
\usepackage{pdfpages}
\usepackage{multicol}

\usetikzlibrary{shapes.geometric, arrows}

\tikzstyle{square} = [rectangle, minimum width=1cm, minimum height = 1cm, text centered, draw = black]

\tikzstyle{diam} = [diamond, minimum width=1cm, minimum height = 1cm, text centered, draw = black]

\tikzstyle{arrow} = [thick, ->, >= stealth]

\usepackage[nameinlink,capitalize]{cleveref}

\usepackage[normalem]{ulem}
\useunder{\uline}{\ul}{}

\usepackage{listings}

\lstdefinestyle{cpp}{
	basicstyle=\small,
	belowcaptionskip=1\baselineskip,
	breaklines=true,
	frame=single,
	language=C++,
	showstringspaces=false,
	basicstyle=\color{black},
	keywordstyle=\bfseries\color{black},
	commentstyle=\itshape\color{white!40!black},
	identifierstyle=\color{blue!60!black},
	stringstyle=\color{green},
	morekeywords={parallel,to},
	escapechar=\& 
}

\usepackage[nottoc,numbib]{tocbibind}

\usepackage{pgf,tikz,pgfplots}
\pgfplotsset{compat=1.15}
\usepackage{mathrsfs}
\usetikzlibrary{arrows}

\pgfplotsset{colormap={cool}{rgb255(0cm)=(255,255,255); rgb255(1cm)=(0,128,255); rgb255(2cm)=(255,0,255)}}

\allowdisplaybreaks[1]

\makeindex

\def\cpp {C\texttt{++} }

\def\dof {DOF }
\def\dofs {DOFs }

   \makeatletter
\let\@fnsymbol\@arabic

\newtheorem{remark}{Remark}
\usepackage{chngcntr}
\counterwithin{remark}{section}

\title{A flexible sparse matrix data format and parallel algorithms for the assembly of sparse matrices in general finite element applications using atomic synchronisation primitives
}
\author{\normalsize{
		Adam Sky\thanks{Corresponding author: Adam Sky, Institute for Structural Mechanics and Dynamics, Technical University Dortmund, August-Schmidt-Str. 8, 44227 Dortmund, Germany, email: adam.sky@tu-dortmund.de}
		\quad , \quad
		César Polindara\thanks{César Polindara, Structural Analysis of Plates and Shells, University of Duisburg-Essen, Universit\"{a}tsstr. 15, 45141 Essen, Germany, email: cesar.polindara-lopez@uni-due.de}
		\quad , \quad
		Ingo Muench\thanks{Ingo Muench, Institute for Structural Mechanics and Dynamics, Technical University Dortmund, August-Schmidt-Str. 8, 44227 Dortmund, Germany, email: ingo.muench@tu-dortmund.de}
		\quad and \quad
		Carolin Birk\thanks{Carolin Birk, Structural Analysis of Plates and Shells, University of Duisburg-Essen, Universit\"{a}tsstr. 15, 45141 Essen, Germany, email: statik-ftw@uni-due.de}
	}
}

\begin{document}

\maketitle

\begin{abstract}
Finite element methods require the composition of the global stiffness matrix from local finite element contributions. The composition process combines the computation of element stiffness matrices and their assembly into the global stiffness matrix, which is commonly sparse. In this paper we focus on the assembly process of the global stiffness matrix and explore different algorithms and their efficiency on shared memory systems using C\texttt{++}. A key aspect of our investigation is the use of atomic synchronization primitives for the derivation of data-race free algorithms and data structures.
Furthermore, we propose a new flexible storage format for sparse matrices and compare its performance with the compressed row storage format using abstract benchmarks based on common characteristics of finite element problems.       
\\
\vspace*{0.25cm}
\end{abstract}

\tableofcontents

\section{Introduction}
The finite element (FE) method approximates solutions to partial differential equations (PDEs) by subdividing a problem domain into a finite number of subdomains, called finite elements. The properties of the PDE in the subdomain are represented by the element's stiffness matrix, the term dates back to the development of the method for elasticity problems but serves to represent the finite element matrices of all kinds of PDEs. In order to find solutions in the domain, the finite element contributions must be computed and assembled into a global stiffness matrix, which is sparse and often symmetric.
This assembly process \cite{JAVAFE,Unterkircher2005} is a time-consuming operation that affects the performance of FE programs, especially when iterative solution schemes are employed. Therefore, the parallelization of the process can lead to vast improvements in performance of the software. 

In modern processor architectures parallelization is available in three forms. The first form is multithreading. A processor contains multiple cores, each of which is capable of handling one or more threads of tasks on its own. Splitting a program into several tasks and spreading them across multiple cores can yield faster computations. In this case, however, synchronization is required to ensure data integrity in case two or more threads attempt to modify the same data simultaneously, a phenomenon known as data race. The second form of parallelization is vectorization, also known as Single Instruction Multiple Data (SIMD). Each processor core is capable of executing a single instruction on a set of multiple data, for which the common requirement is contiguous data alignment. The last form of parallelization is Instructions Per Cycle (IPC). Depending on the instruction and latency, the processor might be able to perform more than one instruction in a single cycle. Latter is not discussed in this work.

Memory efficiency is another major aspect in designing for performance. The main hardware elements that make up the computer memory structure are roughly speaking, in increasing order of response time and capacity, the cache, the RAM and hard drive. Computationally expensive tasks in a FE program like matrix assembly, matrix-vector multiplication or matrix factorization, to name a few, use primarily cache and RAM and therefore are affected by the concept of data locality, often called locality of the reference. This concept refers to the tendency of a computer program to achieve faster access to objects whose addresses are near one another, either in space (spatial locality) or in time (temporal locality). It has been shown \cite{pichel2005performance,Lee:CSD-03-1297,10.1007/11557654_91} that by improving data locality a performance boost in SpMV (Sparse Matrix by dense vector product) algorithms is obtained. We show in this work that the same idea can be applied to matrix assembly algorithms.

Available literature for assembly algorithms mostly considers strategies for dividing the assembly, such that data races do not occur\cite{NovaesDeRezende2000,Unterkircher2005}. Perhaps the most common example is the use of colouring algorithms \cite{CATALYUREK2012576,inproceedings,Natarajan1991,Cecka2011,Jones1993}. Another strategy, known as divide and conquer\cite{Thebault2013,inproceedings,1676280}, suggests splitting the finite element mesh across the processor cores such that bordering elements are assembled separately.
While these approaches allow for a lock-free assembly, they usually require rearrangement of data and the setting of barriers for synchronization. For example, colour graph methods can assemble one colour at a time in parallel and a barrier is set between colours. In other words, an external sequential loop iterates over colours. 

In contrast, this paper focuses on the application of synchronization primitives in order to derive data-race free data structures and algorithms and compares their respective performance. Specifically, the usage of atomics is emphasized. Further, the performance of the new algorithms is compared to a colouring algorithm in order to frame them with respect to the commonly used methods. 

There are many available sparse matrix storage formats for finite element programs \cite{Sander2020} which allow to reduce the storage requirements of the data structure.
Perhaps the most common of which is the compressed storage row (CSR) format. Consequently, we employ the CSR format as our frame of reference for performance measurements. Another common format, specifically for finite element computations, is the block compressed storage row (BCSR) format, which hard-codes specific matrix block sizes into the sparse matrix data structure depending on the problem type. In this paper we introduce a new flexible sparse matrix storage format designated compressed row aligned columns (CRAC). The flexibility of the format stems from its storage of blocks of various sizes, making it independent of the problem type. Consequently, this format integrates well in finite element programs, where the problem to be solved is defined at runtime (dynamically).

\begin{remark}
	The focus of this paper is \textbf{not} the improvement in performance of a \textbf{specific} finite element problem or the definition of a correspondingly matching sparse matrix format, but rather the introduction and a general investigation of the performance of parallel assembly algorithms with respect to the parallelization strategy and the sparse matrix format. Consequently, the computation of each element stiffness matrix is \textbf{not} the main focus of this article and dummy elements are employed. We note that in terms of performance, the computation of the element stiffness matrices is greatly affected by problem type and the employed methodology, compare \cite{Sch2014}.
	\label{re:1}
\end{remark}

This paper is organized as follows: In \cref{sec:data} we present the data structures used to represent the element and global stiffness matrices and introduce a new storage format for sparse matrices. In \cref{sec:ass} we discuss parallelization strategies based on the current \cpp standards and introduce the corresponding pseudocode. Then, in \cref{seq:res} the numerical experiments on a specific architecture, namely, a Windows server, are presented and discussed.
Finally, in \cref{sec:out} we give our conclusions and an outlook for further development.

\section{Data structures for finite element programs}\label{sec:data}
In this section we describe the data structures used for the element stiffness matrix, hereinafter referred to as \lstinline|K|, and the global sparse stiffness matrix, hereinafter referred to as \lstinline|G|.

The assembly process in the finite element method usually results in a global sparse matrix, i.e., a matrix that consists mostly of zeros. As mesh grids containing many elements generate large global matrices and due to the memory limitations of standard computers, many types of storage formats have been developed for the sparse matrix, where the zero coefficients are omitted and only none zero values are stored, see e.g. \cite{10.5555/1855048, Sander2020}.
One of the most widely employed formats is the so called compressed storage row (CSR) and its column-major analogue, the compressed storage column (CSC). Another common format, specifically in the realm of finite elements, is the block compressed storage row (BCSR). The BCSR format stores block matrices of a fixed size and maps to them in the same manner the CSR format maps to values. The usefulness of this features is made apparent by multi-dimensional problems. For example, in plane-stress linear elasticity each node on the finite element mesh has two translational degrees of freedom. If these are numbered contiguously, each node implies a contribution of a dense $\mathbb{R}^{2\times 2}$ matrix to the global stiffness matrix. Consequently, a BCSR format storing blocks of $\mathbb{R}^{2\times2}$ matrices is ideal, see \cite{Sander2020}. 

In the following we explain CSR and introduce a new format, drawing design aspects from both CSR and BCSR. We designate the new format as \textbf{C}ompressed \textbf{R}ow \textbf{A}ligned \textbf{C}olumns (CRAC). Its column-major counterpart is the compressed columns aligned rows (CCAR).
Furthermore, we discuss the generic element and consider possible storage schemes for its \dofs map.

\subsection{A generic element}
In order to assess the efficiency of the assembly algorithms presented in Section 3, the operations associated to the computation of the local stiffness matrix are omitted, see \cref{re:1}. Therefore, solely a dynamically allocated constant local stiffness matrix is defined for each element. Furthermore, since in our experiments we make use of quadrilateral meshes and the number of degrees of freedom (DOFs) per node is constant for each benchmark, we employ a matrix filled with ones as double precision floating-point number in accordance with an element's DOFs. The local stiffness matrix is associated with an array \lstinline|element_dofs_array| containing the \dof indices that map \lstinline|K| into \lstinline|G|. Both data structures, \lstinline|K| and \lstinline|element_dofs_array|, define a generic element.

The \lstinline|element_dofs_array| of a finite element can be stored in two manners. Let us consider, for a four-node quadrilateral element with 4 \dofs per node with a typical but arbitrary sequence of numbers as an example:
\begin{align}
\lstinline|element_dofs_array| := [52 \quad 53\quad 54\quad 55\quad 96\quad 97\quad 98\quad 99\quad 100\quad 101\quad 102\quad 103\quad 48\quad 49\quad 50\quad 51] \, . \notag
\end{align}
If a node has multiple degrees of freedom, the array contains several series of contiguous \dof indices. Furthermore, a node might be neighbouring another, such that their combined \dof indices generate a contiguous sequence. This motivates a new storage scheme, where the advantage of this characteristic is taken into account. Namely, one stores only the first index in a sequence and the number of contiguous consecutive terms, in the example above:
\begin{align}
\lstinline|element_dofs_array| := [52\quad \textcolor{blue}{4}\quad 96\quad \textcolor{blue}{8}\quad 48\quad \textcolor{blue}{4}] \, . \notag
\end{align}
This representation is usually more compact, allowing for a reduction in data storage. Additionally, this representation becomes advantageous in Section 3.3, where vectorization is discussed.

\subsection{The CSR and CRAC formats}
The CSR format stores a matrix in three arrays. The first is the row pointers array, containing the indices for the start of each row in the columns and values arrays.
Correspondingly, the end of a row is given by the following integer in the row pointers array. As such, the row pointers array requires one row more than the matrix in order to store the end of the last row. The columns and values arrays have the same size. The columns array maps values to their respective column index as both share the same position in their respective arrays. 
 
The BCSR format is used primarily in FE-software, since for two or more \dof per node the \dofs belonging to one node yield a contiguous sequence of fixed size, and their stiffness terms can be stored as a block matrix in the global sparse matrix, i.e., a small sub matrix.    
The CRAC format extends the idea presented in the element \dof indices array and relaxes the assumptions of the BCSR format by introducing row blocks that vary in length. 

\begin{remark}
		Note that unlike in BCSR, the CRAC format does not actually store block matrices, but rather the start and end of a column block in the values array. Consequently, the data in CRAC remains contiguous and compact, whereas BCSR may lead to padding if the sub matrix is not a multiple of 16 byte memory blocks.
\end{remark}

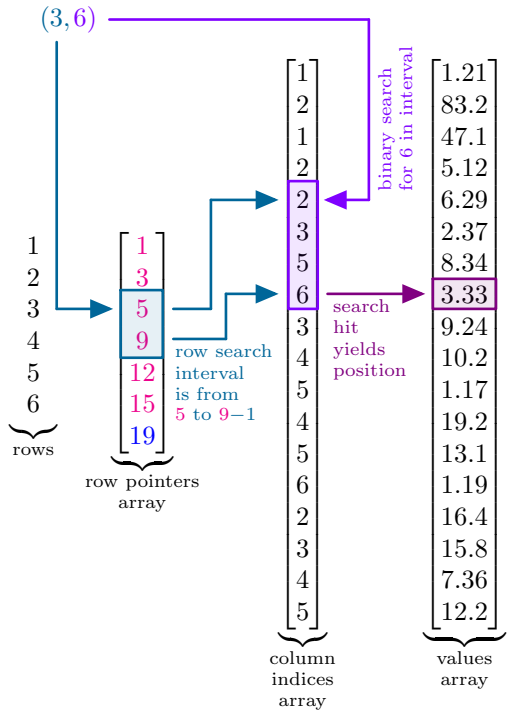
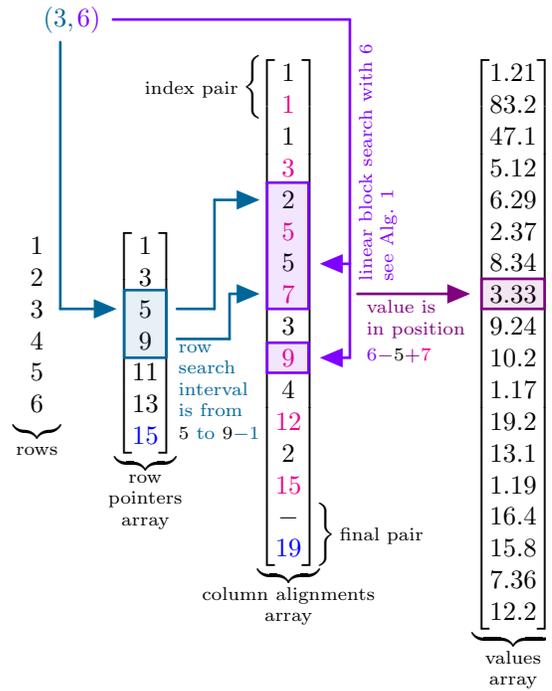
\begin{figure}
	\begin{subfigure}{0.48\linewidth}
		\centering
		\definecolor{yqqqyq}{rgb}{0.5019607843137255,0.,0.5019607843137255}
\begin{tikzpicture}[scale = 0.2][line cap=round,line join=round,>=triangle 45,x=1.0cm,y=1.0cm]
\clip(5,15) rectangle (45,45);
\draw [line width=1.5pt] (10.,45.)-- (10.,15.);
\draw [line width=1.5pt] (10.,15.)-- (40.,15.);
\draw [line width=1.5pt] (40.,15.)-- (40.,45.);
\draw [line width=1.5pt] (10.,45.)-- (40.,45.);
\draw [line width=1.5pt] (15.,45.)-- (15.,15.);
\draw [line width=1.5pt] (20.,45.)-- (20.,15.);
\draw [line width=1.5pt] (25.,45.)-- (25.,15.);
\draw [line width=1.5pt] (30.,45.)-- (30.,15.);
\draw [line width=1.5pt] (35.,45.)-- (35.,15.);
\draw [line width=1.5pt] (10.,40.)-- (40.,40.);
\draw [line width=1.5pt] (10.,35.)-- (40.,35.);
\draw [line width=1.5pt] (10.,30.)-- (40.,30.);
\draw [line width=1.5pt] (10.,25.)-- (40.,25.);
\draw [line width=1.5pt] (10.,20.)-- (40.,20.);
\draw (10.5,43.5) node[anchor=north west] {$1.21$};
\draw (21,43.5) node[anchor=north west] {$0$};
\draw (15.5,43.5) node[anchor=north west] {$83.2$};
\draw (26,43.5) node[anchor=north west] {$0$};
\draw (31,43.5) node[anchor=north west] {$0$};
\draw (36,43.5) node[anchor=north west] {$0$};
\draw (10.5,38.5) node[anchor=north west] {$47.1$};
\draw (21,38.5) node[anchor=north west] {$0$};
\draw (15.5,38.5) node[anchor=north west] {$5.12$};
\draw (26,38.5) node[anchor=north west] {$0$};
\draw (31,38.5) node[anchor=north west] {$0$};
\draw (36,38.5) node[anchor=north west] {$0$};
\draw (11,33.5) node[anchor=north west] {$0$};
\draw (15.5,33.5) node[anchor=north west] {$6.29$};
\draw (26,33.5) node[anchor=north west] {$0$};
\draw (20.5,33.5) node[anchor=north west] {$2.37$};
\draw (35.5,33.5) node[anchor=north west] {$3.33$};
\draw (30.5,33.5) node[anchor=north west] {$8.34$};
\draw [line width=1.pt,color=yqqqyq] (35.7,33.5) -- (39.5,33.5) -- (39.5,31.4) -- (35.7,31.4) --(35.7,33.5);
\fill[line width=2.pt,color=yqqqyq,fill=yqqqyq,fill opacity=0.10000000149011612] (35.7,33.5) -- (39.5,33.5) -- (39.5,31.4) -- (35.7,31.4) -- cycle;
\draw (41,24) node[anchor=north west,rotate=90,color=yqqqyq] {$_{\text{Coefficient of our example}}$};
\draw [line width=1.pt,color=yqqqyq] (41,24) -- (43,24) -- (43,22) -- (41,22) --(41,24);
\fill[line width=2.pt,color=yqqqyq,fill=yqqqyq,fill opacity=0.10000000149011612](41,24) -- (43,24) -- (43,22) -- (41,22) -- cycle;
\draw (11,28.5) node[anchor=north west] {$0$};
\draw (20.5,28.5) node[anchor=north west] {$9.24$};
\draw (16,28.5) node[anchor=north west] {$0$};
\draw (25.5,28.5) node[anchor=north west] {$10.2$};
\draw (30.5,28.5) node[anchor=north west] {$1.17$};
\draw (36,28.5) node[anchor=north west] {$0$};
\draw (11,23.5) node[anchor=north west] {$0$};
\draw (21,23.5) node[anchor=north west] {$0$};
\draw (16,23.5) node[anchor=north west] {$0$};
\draw (25.5,23.5) node[anchor=north west] {$19.2$};
\draw (30.5,23.5) node[anchor=north west] {$13.1$};
\draw (35.5,23.5) node[anchor=north west] {$1.19$};
\draw (11,18.5) node[anchor=north west] {$0$};
\draw (20.5,18.5) node[anchor=north west] {$15.8$};
\draw (15.5,18.5) node[anchor=north west] {$16.4$};
\draw (25.5,18.5) node[anchor=north west] {$7.36$};
\draw (30.5,18.5) node[anchor=north west] {$12.2$};
\draw (36,18.5) node[anchor=north west] {$0$};
\end{tikzpicture}
		\caption{}
		\label{fig:mat}
	\end{subfigure}
	\begin{subfigure}{0.48\linewidth}
		\centering
		\definecolor{qqqqff}{rgb}{0.407, 0.776, 0.874}
\begin{tikzpicture}[scale = 0.2][line cap=round,line join=round,>=triangle 45,x=1.0cm,y=1.0cm]
\clip(50,15) rectangle (80,45);
\fill[line width=2.pt,color=qqqqff,fill=qqqqff,fill opacity=0.5] (50.,45.) -- (60.,45.) -- (60.,35.) -- (50.,35.) -- cycle;
\fill[line width=2.pt,color=qqqqff,fill=qqqqff,fill opacity=0.5] (55.,35.) -- (65.,35.) -- (65.,30.) -- (55.,30.) -- cycle;
\fill[line width=2.pt,color=qqqqff,fill=qqqqff,fill opacity=0.5] (70.,35.) -- (80.,35.) -- (80.,30.) -- (70.,30.) -- cycle;
\fill[line width=2.pt,color=qqqqff,fill=qqqqff,fill opacity=0.5] (60.,30.) -- (75.,30.) -- (75.,25.) -- (60.,25.) -- cycle;
\fill[line width=2.pt,color=qqqqff,fill=qqqqff,fill opacity=0.5] (65.,25.) -- (80.,25.) -- (80.,20.) -- (65.,20.) -- cycle;
\fill[line width=2.pt,color=qqqqff,fill=qqqqff,fill opacity=0.5] (55.,20.) -- (75.,20.) -- (75.,15.) -- (55.,15.) -- cycle;
\draw [line width=1.5pt] (50.,45.)-- (50.,15.);
\draw [line width=1.5pt] (50.,15.)-- (80.,15.);
\draw [line width=1.5pt] (80.,15.)-- (80.,45.);
\draw [line width=1.5pt] (80.,45.)-- (50.,45.);
\draw [line width=1.5pt] (55.,15.)-- (55.,45.);
\draw [line width=1.5pt] (60.,45.)-- (60.,15.);
\draw [line width=1.5pt] (65.,45.)-- (65.,15.);
\draw [line width=1.5pt] (70.,45.)-- (70.,15.);
\draw [line width=1.5pt] (75.,45.)-- (75.,15.);
\draw [line width=1.5pt] (50.,40.)-- (80.,40.);
\draw [line width=1.5pt] (50.,35.)-- (80.,35.);
\draw [line width=1.5pt] (50.,30.)-- (80.,30.);
\draw [line width=1.5pt] (50.,25.)-- (80.,25.);
\draw [line width=1.5pt] (50.,20.)-- (80.,20.);
\draw (51,43.5) node[anchor=north west,color=magenta] {$1$};
\draw (56,43.5) node[anchor=north west] {$2$};
\draw (51,38.5) node[anchor=north west,color=magenta] {$3$};
\draw (56,38.5) node[anchor=north west] {$4$};
\draw (56,33.5) node[anchor=north west,color=magenta] {$5$};
\draw (61,33.5) node[anchor=north west] {$6$};
\draw (71,33.5) node[anchor=north west,color=magenta] {$7$};
\draw (76,33.5) node[anchor=north west] {$8$};
\draw (61,28.5) node[anchor=north west,color=magenta] {$9$};
\draw (66,28.5) node[anchor=north west] {$10$};
\draw (71,28.5) node[anchor=north west] {$11$};
\draw (66,23.5) node[anchor=north west,color=magenta] {$12$};
\draw (71,23.5) node[anchor=north west] {$13$};
\draw (76,23.5) node[anchor=north west] {$14$};
\draw (56,18.5) node[anchor=north west,color=magenta] {$15$};
\draw (61,18.5) node[anchor=north west] {$16$};
\draw (66,18.5) node[anchor=north west] {$17$};
\draw (71,18.5) node[anchor=north west] {$18$};
\end{tikzpicture}
		\caption{}
		\label{fig:pat}
	\end{subfigure}
	\begin{subfigure}{0.48\linewidth}
		\centering
		\definecolor{qqwwzz}{rgb}{0.,0.4,0.6}
		\definecolor{xfqqff}{rgb}{0.4980392156862745,0.,1.}
		\definecolor{yqqqyq}{rgb}{0.5019607843137255,0.,0.5019607843137255}
	\begin{tikzpicture}[line cap=round,line join=round,>=triangle 45,x=1.0cm,y=1.0cm]
					\clip(0,-0.5) rectangle (8,9.2);
					\draw (1.4,9) node[anchor=north west,color=qqwwzz] {$(3,$};
					\draw (1.85,9) node[anchor=north west,color=xfqqff] {$6)$};
					\draw (1,6)  node[anchor=north west] {$\underbrace{\begin{matrix}
								1 \\ 2 \\ 3 \\ 4 \\ 5 \\ 6 \\ 
						\end{matrix}}_{\text{rows}}$};
					\draw (2,6)  node[anchor=north west] {$\underbrace{\begin{bmatrix}
								\textcolor{magenta}{1} \\ \textcolor{magenta}{3} \\ \textcolor{magenta}{5} \\ \textcolor{magenta}{9} \\ \textcolor{magenta}{12} \\ \textcolor{magenta}{15} \\ 
								\textcolor{blue}{19}
						\end{bmatrix}}_{\parbox{4.35em}{\centering \scriptsize row pointers \\  array}}$};
					\draw (4,8.3)  node[anchor=north west] {$\underbrace{\begin{bmatrix}
								1 \\ 2 \\ 1 \\ 2 \\ 2 \\ 3 \\ 5 \\ 6 \\ 3 \\ 4 \\ 5
								\\ 4 \\ 5 \\ 6 \\ 2 \\ 3 \\ 4 \\ 5
						\end{bmatrix}}_{\parbox{5.1em}{\centering \scriptsize column indices \\  array}}$};
					\draw (6.5,8.3)  node[anchor=north west] {$\underbrace{\begin{bmatrix}
								1.21 \\ 83.2 \\ 47.1 \\ 5.12 \\ 6.29 \\ 2.37 \\ 8.34 \\ 3.33 \\ 9.24 
								\\ 10.2 \\ 1.17 \\ 19.2 \\13.1 \\ 1.19 \\ 16.4 \\ 15.8 \\ 7.36 \\ 12.2
						\end{bmatrix}}_{\parbox{3em}{\centering \scriptsize values \\  array}}$};
					\draw [line width=1.pt,color=qqwwzz] (1.75,8.4) -- (1.75,4.85);
					\draw [->,line width=1.pt,color=qqwwzz] (1.75,4.85) -- (2.5,4.85);
					\draw [line width=1.pt,color=qqwwzz] (2.6,5.1) -- (3.16,5.1) -- (3.16,4.2) -- (2.6,4.2) -- (2.6,5.1);
					\fill[line width=2.pt,color=qqwwzz,fill=qqwwzz,fill opacity=0.10000000149011612] (2.6,5.1) -- (3.16,5.1) -- (3.16,4.2) -- (2.6,4.2) -- cycle;
					\draw (3.2,4.5) node[anchor=north west,color=qqwwzz] {\parbox{5em}{ \scriptsize row search \\  interval \\ is from} };
					\draw (3.2,3.65) node[anchor=north west,color=qqwwzz] {$_{\textcolor{magenta}{5} \text{ to } \textcolor{magenta}{9}-1}$ };
					\draw [line width=1.pt,color=qqwwzz] (3.3,4.85) -- (3.8,4.85);
					\draw [line width=1.pt,color=qqwwzz] (3.8,4.85) -- (3.8,6.3);
					\draw [->,line width=1.pt,color=qqwwzz] (3.8,6.3) -- (4.7,6.3);
					\draw [line width=1.pt,color=qqwwzz] (3.3,4.45) -- (4,4.45);
					\draw [line width=1.pt,color=qqwwzz] (4,4.45) -- (4,5.05);
					\draw [->,line width=1.pt,color=qqwwzz] (4,5.05) -- (4.7,5.05);
					\draw [line width=1.pt,color=xfqqff] (4.83,6.55) -- (5.22,6.55) -- (5.22,4.85) -- (4.83,4.85) -- (4.83,6.55);
					\fill[line width=2.pt,color=xfqqff,fill=xfqqff,fill opacity=0.10000000149011612] (4.83,6.55) -- (5.22,6.55) -- (5.22,4.85) -- (4.83,4.85) -- cycle;
					\draw [line width=1.pt,color=xfqqff] (2.45,8.7) -- (5.9,8.7);
					\draw [line width=1.pt,color=xfqqff] (5.9,8.7) -- (5.9,6.3);
					\draw (5.9,6.3) node[anchor=north west,rotate=90,color=xfqqff] {$\parbox{6em}{\scriptsize binary search \\ for $6$ in interval}$};
					\draw [->,line width=1.pt,color=xfqqff] (5.9,6.3) -- (5.3,6.3);
					\draw [->,line width=1.pt,color=yqqqyq] (5.35,5.05) -- (6.65,5.05);
					\draw (5.3,5.1) node[anchor=north west,color=yqqqyq] {$\parbox{3em}{ \scriptsize search hit \\  yields position}$};
					\draw [line width=1.pt,color=yqqqyq] (6.75,5.25) -- (7.59,5.25) -- (7.59,4.85) -- (6.75,4.85) -- (6.75,5.25);
					\fill[line width=2.pt,color=yqqqyq,fill=yqqqyq,fill opacity=0.10000000149011612] (6.75,5.25) -- (7.59,5.25) -- (7.59,4.85) -- (6.75,4.85) -- cycle;
	\end{tikzpicture}
		\caption{}
		\label{fig:csr}
	\end{subfigure}
\begin{subfigure}{0.48\linewidth}
	\centering
	\definecolor{qqwwzz}{rgb}{0.,0.4,0.6}
	\definecolor{xfqqff}{rgb}{0.4980392156862745,0.,1.}
	\definecolor{yqqqyq}{rgb}{0.5019607843137255,0.,0.5019607843137255}
	\begin{tikzpicture}[line cap=round,line join=round,>=triangle 45,x=1.0cm,y=1.0cm]
		\clip(0,-0.5) rectangle (8.2,9.2);
		\draw (1.4,9) node[anchor=north west,color=qqwwzz] {$(3,$};
		\draw (1.85,9) node[anchor=north west,color=xfqqff] {$6)$};
		\draw (1,6)  node[anchor=north west] {$\underbrace{\begin{matrix}
					1 \\ 2 \\ 3 \\ 4 \\ 5 \\ 6 \\ 
			\end{matrix}}_{\text{rows}}$};
		\draw (2,6)  node[anchor=north west] {$\underbrace{\begin{bmatrix}
					1 \\ 3 \\ 5 \\ 9 \\ 11 \\ 13 \\ \textcolor{blue}{15}
			\end{bmatrix}}_{\parbox{4.3em}{\centering \scriptsize row pointers \\ array}}$};
		\draw (2.75,8.35)  node[anchor=north west] {$\text{\scriptsize index pair} \left \{  \begin{matrix} \phantom{1} \\ \phantom{1} \end{matrix} \right . $};
		\draw (4.8,2.4)  node[anchor=north west] {$ \left .  \begin{matrix} \phantom{1} \\ \phantom{1} \end{matrix} \right \} \text{\scriptsize final pair}  $};
		\draw (3.5,8.3)  node[anchor=north west] {$\underbrace{\begin{bmatrix}
					1 \\ \textcolor{magenta}{1} \\ 
					1 \\ \textcolor{magenta}{3} \\
					2 \\ \textcolor{magenta}{5} \\ 
					5 \\ \textcolor{magenta}{7} \\ 
					3 \\ \textcolor{magenta}{9} \\ 
					4 \\ \textcolor{magenta}{12} \\
					2 \\ \textcolor{magenta}{15} \\ 
					- \\ \textcolor{blue}{19} 
			\end{bmatrix}}_{\parbox{6.6em}{\centering \scriptsize column alignments \\ array}}$};
		\draw (7.1,8.3)  node[anchor=north west] {$\underbrace{\begin{bmatrix}
					1.21 \\ 83.2 \\ 47.1 \\ 5.12 \\ 6.29 \\ 2.37 \\ 8.34 \\ 3.33 \\ 9.24 
					\\ 10.2 \\ 1.17 \\ 19.2 \\13.1 \\ 1.19 \\ 16.4 \\ 15.8 \\ 7.36 \\ 12.2
			\end{bmatrix}}_{\parbox{3em}{\centering \scriptsize values \\  array}}$};
		\draw [line width=1.pt,color=qqwwzz] (1.75,8.4) -- (1.75,4.85);
		\draw [->,line width=1.pt,color=qqwwzz] (1.75,4.85) -- (2.5,4.85);
		\draw [line width=1.pt,color=qqwwzz] (2.6,5.1) -- (3.16,5.1) -- (3.16,4.2) -- (2.6,4.2) -- (2.6,5.1);
		\fill[line width=2.pt,color=qqwwzz,fill=qqwwzz,fill opacity=0.10000000149011612] (2.6,5.1) -- (3.16,5.1) -- (3.16,4.2) -- (2.6,4.2) -- cycle;
		\draw (3.2,4.5) node[anchor=north west,color=qqwwzz] {\parbox{3em}{ \scriptsize row \\ search \\  interval \\ is from} };
		\draw (3.2, 3.4) node[anchor=north west,color=qqwwzz] {
		$_{\textcolor{black}{5} \text{ to } \textcolor{black}{9}-1}$};
		\draw [line width=1.pt,color=qqwwzz] (3.3,4.85) -- (3.8,4.85);
		\draw [line width=1.pt,color=qqwwzz] (3.8,4.85) -- (3.8,6.3);
		\draw [->,line width=1.pt,color=qqwwzz] (3.8,6.3) -- (4.4,6.3);
		\draw [line width=1.pt,color=qqwwzz] (3.3,4.45) -- (4,4.45);
		\draw [line width=1.pt,color=qqwwzz] (4,4.45) -- (4,5.05);
		\draw [->,line width=1.pt,color=qqwwzz] (4,5.05) -- (4.4,5.05);
		\draw [line width=1.pt,color=xfqqff] (4.5,6.53) -- (5.05,6.53) -- (5.05,4.85) -- (4.5,4.85) -- (4.5,6.53);
		\draw [line width=1.pt,color=xfqqff] (4.5,4.) -- (5.05,4.) -- (5.05,4.4) -- (4.5,4.4) -- (4.5,4.);
		\fill[line width=2.pt,color=xfqqff,fill=xfqqff,fill opacity=0.10000000149011612] (4.5,6.53) -- (5.05,6.53) -- (5.05,4.85) -- (4.5,4.85) -- cycle;
		\fill[line width=2.pt,color=xfqqff,fill=xfqqff,fill opacity=0.10000000149011612] (4.5,4.) -- (5.05,4.) -- (5.05,4.4) -- (4.5,4.4) -- cycle;
		\draw [line width=1.pt,color=xfqqff] (2.45,8.7) -- (5.6,8.7);
		\draw [line width=1.pt,color=xfqqff] (5.6,8.7) -- (5.6,4.2);
		\draw [->,line width=1.pt,color=xfqqff] (5.6,5.45) -- (5.2,5.45);
		\draw [->,line width=1.pt,color=xfqqff] (5.6,4.2) -- (5.2,4.2);
		\draw (5.6,5.1) node[anchor=north west,rotate=90,color=xfqqff] {${\parbox{10.em}{ \scriptsize linear block search with 6 \\ see Alg. 1}}$};
		\draw (5.7,5.1) node[anchor=north west,color=yqqqyq] {$_{\text{value is}}$};
		\draw (5.7,4.8) node[anchor=north west,color=yqqqyq] {$_{\text{in position}}$};
		\draw (5.7,4.45) node[anchor=north west,color=yqqqyq] {$_{ \textcolor{xfqqff}{6}- \textcolor{black}{5} +  \textcolor{magenta}{7} }$};
		\draw [line width=1.pt,color=yqqqyq] (7.35,4.85) -- (8.19,4.85) -- (8.19,5.25) -- (7.35,5.25) -- (7.35,4.85);
		\fill[line width=2.pt,color=yqqqyq,fill=yqqqyq,fill opacity=0.10000000149011612] (7.35,4.85) -- (8.19,4.85) -- (8.19,5.25) -- (7.35,5.25) -- cycle;
		\draw [->,line width=1.pt,color=yqqqyq] (5.7,5.05) -- (7.2,5.05);
	\end{tikzpicture}
	\caption{}
	\label{fig:crac}
\end{subfigure}
\caption{Standard matrix representation (a) and the corresponding sparsity pattern with indices (b). Search for coefficient $(3,6)$ in CSR (c) format and CRAC (d) format representations.}
\label{fig:formats}
\end{figure}

Instead of storing the indices of all columns, the column alignments array stores the first index of consecutive contiguous coefficients and its corresponding position in the values array. I.e., the column alignments array stores pairs of integers, see example in \cref{fig:formats}. As such, the column alignments array requires at least two entries per row and one final pair to mark the end of the values array. 

We compare the access pattern of the CSR and CRAC formats in retrieving the coefficient $(3,6)$ as an example, see also \cref{fig:formats}.

\begin{lstlisting}[caption =  Linear block search for the CRAC format]
	i = row_start_position //start of the row in the column alignments array
	c = column_index //the index of the searched column
	while i < row_end
	  next_column_index = column_alignments_array(i + 2)
	  if next_column_index > c //compare with next column index
	    break //break the loop if larger
	  i = i + 2 
	column_start_index = column_alignments_array(i) //first column block index 
	values_array_start_position = column_alignments_array(i + 1) //column block start in values array
	values_array_end_position = column_alignments_array(i + 3)
	colum_end_index = column_start_index + (values_array_end_position - values_array_start_position)
	//compute position if in block
	if c >= column_start_index and c < colum_end_index
	  return values_array_start_position + (c - column_start_index)
\end{lstlisting}

In the CSR format, ones enters the row pointers array at positions $3$ and $3+1 = 4$ and retrieves the indices $5$ and $9$. The column indices of that row are therefore stored in the column indices array between the positions $5$ and $9-1 = 8$. Searching the column indices array between $5$ and $8$ for the number $6$ yields the position $8$. Consequently, the value of the coefficient $(3,6)$ is found at the position $8$ of the values array, reading $3.33$.

In the analogous procedure in the CRAC format one enters the row pointers array at positions $3$ and $3+1 = 4$, retrieving the indices $5$ and $9$. Therefore, the data regarding the row is contained in the column alignments array between the positions $5$ and $9-1 = 8$. One searches the column alignments array between the positions 5 and 8 for an index larger than 6 using jumps of 2. In other words man compares only the starting indices of column blocks with the column number 6. If a larger index is found, the previous block may contain the searched column number. Else the last column block of the row is fetched, as is the case in our example. Consequently, one finds the indices $5$ and $7$. The first index represents the start of the relevant column block and the second index is its position in the values array. The starting position of the next column block in the values array is given at $8+2 = 10$ of the column alignments array and contains the integer $9$. Consequently, the column block goes from $5$ to $5 + 9 - 7 - 1 = 6$, i.e. the length of the column block is $9-7 = 2$. As the searched column index $6$ is in the column blocks range ($ 5 \leq 6 \leq 6$), the position of its corresponding value can be calculated by subtracting the difference in column indices and adding it to the position of the values array $6 - 5 + 7  = 8$. Thus, the value is found at position $8$ of the values array and reads $3.33$. 

When accessing coefficients in the CSR format, we note the necessity of searching the row for the corresponding column. Depending on the amount of coefficients in a row, we employ either a linear or a binary search. Our investigations show the binary search to outperform the linear search only for rows with more than 30 coefficients. 
While the search algorithm in the CRAC format adds a small layer of complexity, it in-fact allows for faster access speeds, as no search is needed inside a column block, but rather only between  column blocks. Finally, Our investigation shows that using a linear search between column blocks of the CRAC format yields the highest access speed due to the limited number of column blocks per row. This is a direct result of the alignment expected in FE-programs.
     
\section{Assembly algorithms}\label{sec:ass}
In this section we present several assembly algorithms. The initialization of the sparse matrix including its sparsity pattern is performed only once in the preprocessing stage, i.e. the rows- and columns arrays are initialized with the corresponding indices (see \cref{fig:pat}) and the values array is sized to the number of none zero entries (NNZ), where each entry is set to zero. Between assemblies, the entries of the values array are reset to zero. 

\begin{remark}
		Both the initialization and the reset of the global sparse matrix are not measured in our benchmarks, since we focus on the assembly procedure.
\end{remark}

The assembly is done via loops, where \lstinline|K| is assembled into \lstinline|G|.
Improvements to the assembly algorithms will be introduced progressively following the ideas of parallelization and vectorization. 
 
\subsection{Sequential assembly method}
We introduce the sequential method of assembly as a basis for performance comparison. The assembly is done in three loops, the outer loop over the elements array, and two inner loops over rows and columns.
The sequential method takes advantage of neither multi-threading nor vectorization. Hence, it will be used to measure the improvement in performance obtained with the other methods presented in this work. Some compilers try to vectorize simple loops, for which they require data alignment of the coefficients. Since the degrees of freedom of one element are not guaranteed to be aligned in the sparse matrix, no auto-vectorization can take place. 
\begin{lstlisting}[caption = Sequential assembly]
	for e = 1 to number of elements
	  K = stiffness matrix of element e 
	  D = element_dofs_array of element e
	  for r = 1 to number of rows in K
	    for c = 1 to number of columns in K
	      do G(D(r),D(c)) = G(D(r),D(c)) + K(r,c)
\end{lstlisting}

\subsection{Parallelization} 
There are several methods of synchronizing multi-threaded applications in order to prevent data races, e.g. a standard mutex, a colour graph \cite{CATALYUREK2012576} or D\&C strategies (Divide and conquer) \cite{inproceedings}. The C\texttt{++}11 standard introduced new types of synchronization primitives called atomics \cite{CIA}. Atomic primitives are accessed and modified using atomic instructions. These are special compiler instructions, where an action can be seen as either complete or not yet started by all threads; no intermediate state can be observed. Consequently, the modification of atomic variables is inherently thread safe and lock-free. While several atomic instructions are available, for our purposes we consider the atomic \lstinline|exchange|, \lstinline|compare_exchange| and \lstinline|fetch_add| operations. The first instruction simply replaces the value of a target parameter with a new value atomically. The second instruction reads the current value of the parameter and compares it to a different given value. If the values match, the instruction will replace the old value with a third given value and return \lstinline|true| as an indicator of a successful exchange. If the values do not match, no replacement takes place and the instruction returns \lstinline|false| for a failed exchange. The third instruction adds a value to the target atomically. The first two instructions build the foundation for the design of more complex data-race free data-structures and methods \cite{CIA}. 

The C\texttt{++}20 standard introduces full-blown atomic double precision floating-point numbers by implementing the \lstinline|fetch_add| instruction for them. The instruction was not available for non-integers in previous C\texttt{++} standards. These types can be modified concurrently in a lock free manner.

With this first tool at hand, the simplest form of parallel assembly using atomics consists in modifying all the coefficients of the sparse matrix atomically, where the first loop over the elements is now done in parallel\footnote{the \lstinline|parallel for| is done using the parallel \lstinline|std::for_each(std::execution::par_unseq,..)| method introduced in C\texttt{++}17.}.
\begin{lstlisting}[caption=Parallel assembly using atomic addition]
parallel for e = 1 to number of elements //parallel loop over elements
  K = stiffness matrix of element e 
  D = element_dofs_array of element e
  for r = 1 to number of rows in K
    for c = 1 to number of columns in K
      do fetch_add(G(D(r),D(c)), K(r,c)) //modify each coefficient atomically 
\end{lstlisting}

While the parallel atomic method guarantees data-race free assembly, it is not the only solution. Another simple method can be achieved by using an array of standard \lstinline|mutex| and locking each row before modification. This approach would require the matrix to contain another array with the size $80 \,n_{\textrm{rows}}$ ($\,n_{\textrm{rows}}$ is the number of rows in the sparse matrix), since each standard \lstinline|mutex| weighs 80 bytes. For large matrices this can be expensive and in fact, unnecessary. Assuming the time spent working on each row is negligible, a \lstinline|spin_mutex| could be used instead. The \lstinline|spin_mutex| is based on the \lstinline|atomic_flag| type, which can only have two states: locked or unlocked; locking and unlocking is done atomically. A thread will keep trying to lock a \lstinline|spin_mutex| in a loop until it succeeds, effectively causing the thread to wait. The \lstinline|spin_mutex| weighs only 4 bytes allowing to considerably reduce the size of the mutex array to $4 \,n_{\textrm{rows}}$. This is already a significant improvement. Note however, that using another array of locks will always increase the size of the sparse matrix data structure. Instead, by considering the functionality of the \lstinline|spin_mutex| it is also possible to drop the new mutex array altogether. As a \lstinline|spin_mutex| only differs between two states, one can achieve the same effect by considering negative and positive integers
\begin{lstlisting}[caption=Locking an integer]
v = abs(t) //get the absolute value of target t
while(t.compare_exchange(v,-v) == false) //attempt to set negative value
\end{lstlisting}	
\begin{lstlisting}[caption=Unlocking an integer]
v = abs(t) //get the absolute value of target t
t.exchange(v) //replace with absolute value
\end{lstlisting}
Locking the index is achieved by setting its sign to minus. Unlocking is done by replacing the index with its absolute value. Now the \lstinline|spin_int| type can be used to construct the row pointers array of both the CSR and CRAC formats \textbf{without the need} for any additional arrays for synchronization. In other words, the row pointers array now serves both its orginal purpose and as a synchronization mechanism to allow for mutual exclusion of rows. Whenever the index is required for array access, the absolute value is returned.

The parallel assembly with a row mutex now takes the form\footnote{the \lstinline|parallel for| is now done using the parallel \lstinline|std::for_each(std::execution::par,..)| in order to account for locks.}
\begin{lstlisting}[caption=Parallel assembly using a spin mutex]
parallel for e = 1 to number of elements
  K = stiffness matrix of element e 
  D = element_dofs_array of element e
  for r = 1 to number of rows in K
    lock(D(r)) //lock the row of the global matrix
    for c = 1 to number of columns in K
      do G(D(r),D(c)) = G(D(r),D(c)) + K(r,c)
    unlock(D(r)) //unlock the row
\end{lstlisting}

\subsection{Vectorization}
By employing the second format of the \lstinline|element_dofs_array| for the elements, the next form of parallelization, namely vectorization using SIMD instructions, can be used. The SIMD-AVX instruction set is already widely available on modern processors. The instruction set allows adding multiple double precision floating-point numbers to multiple double precision floating-point numbers in a single instruction. The number of double precision floating-point numbers added in a single operation depends on the specific AVX instruction set.

The main requirement of SIMD instructions is memory alignment. Since the DOFs are stored in a row major fashion in a single array, any contiguous sequence of columns is aligned in memory.
Using the second format of the dummy element we can now iterate over aligned column sets and modify the matrix with vector instructions, where the iteration over the element \dofs array is now in jumps of 2, such that two values are extracted from the array at each iteration, the first one being the column index and second one being the size of the contiguous slice.  
By locking each row for modification using the \lstinline|spin_int| version of the sparse matrix storage format, we write the result via SIMD directly to the values array
\begin{lstlisting}[caption=Vectorized parallel assembly with a \lstinline|spin_int| sparse matrix data structure]
parallel for e = 1 to number of elements 
  K = stiffness matrix of element e  
  D = element_dofs_array of element e //D is now the array of dofs and alignment
  for r = 1 to number of rows in K
    lock(D(r)) //lock the row of the global matrix
    c = 1 //for iteration over columns of K
    a = 1 //for iteration over the alignment in D
    while(c < number of columns in K) do
      //modify aligned coefficients via simd instructions
      G(D(r),D(a):D(a) + D(a+1)) = G(D(r),D(a):D(a) + D(a+1)) + K(r,c : c + D(a+1))
      c = c + D(a+1)
      a = a + 2 //start of the next slice in D
    unlock(D(r)) //unlock the row    
\end{lstlisting}
where the $:$ operator indicates the fetching of a slice.

In order to frame our methods with respect to commonly used algorithms, we consider the assembly algorithm with parallelization using colouring while employing the vectorization method

\vspace{3.cm}
\begin{lstlisting}[caption=Parallel assembly using colouring and vectorization]
for colour = 1 to number of colours
  parallel for e = 1 to number of elements in colour
    K = stiffness matrix of element e  
    D = element_dofs_array of element e //D is now the array of dofs and alignment
    for r = 1 to number of rows in K
      c = 1 //for iteration over columns of K 
      a = 1 //for iteration over the alignment in D
      while(c < number of columns in K) do
        //modify aligned coefficients via simd instructions
        G(D(r),D(a):D(a) + D(a+1)) = G(D(r),D(a):D(a) + D(a+1)) + K(r,c : c + D(a+1))
        c = c + D(a+1)
        a = a + 2 //start of the next slice in D
\end{lstlisting}
where the greedy colouring algorithm \cite{10.1007/11823285_61} has been used to generate colours.
Since elements of the same colour do not share any neighbouring nodes, each colour can be assembled in parallel without locks.
  
\section{Performance benchmarks}\label{seq:res}
In this work we quantify the efficiency of an assembly algorithm by the speed factor $c$. To obtain the value, a given assembly algorithm is called $n$ times and its execution time $t_i$ is measured at every call. Thus, the average execution time for a given method is computed as 
	\begin{align}
		t_\text{avg} = \dfrac{1}{n} \sum_i^n t_i \, . 
	\end{align}
The average speed factor for every method except the sequential one is defined accordingly  
	\begin{align}
		c = &\dfrac{t_\text{avg}^\text{seq}}{t_\text{avg}^\text{method}} \, .
\end{align}

Further, we define multiple benchmarks based on the abstract characteristics of finite element meshes for various variational problems. 
The general approach for achieving solution convergence in finite element applications is by mesh refinement. Refinement of the element size, the so called h-refinement, is reflected by the h-benchmark, where we successively reduce the element size. Another strategy is to increase the polynomial power of the base functions, leading to p-refinement and the corresponding p-benchmark, where the amount of nodes on each element is increased. We conduct the p-benchmark by increasing the element power from 1 to 8.
The last type of benchmark is the dimensional benchmark (d-benchmark). This benchmark relates the dimensionality of the solution variable to the amount of degrees of freedom on each node of the finite element mesh. For example, in plane-stress linear elasticity each node has two translational degrees of freedom reflecting the displacement field $u \in \mathbb{R}^d ,\, d = 2$. However, for many types of problems the amount of \dof per node may be much larger. Common examples in mechanics are shell elements \cite{Wallner2020} and generalized continuum theories \cite{Sky}. The d-benchmark is done by increasing the amount of DOF per node from 1 to 8.
In order to characterize the h- and p-benchmarks for problems involving vectorial solutions ($u \in \mathbb{R}^d, \, d > 1$) we also test the methods for the case $d = 4$.

\begin{remark}
		The h-, p- and d-benchmarks all correspond with an increase in size of the global stiffness matrix. However, they lead to different sparsity patterns and differing lock-unlock patterns. For example, in the h-benchmark the size of the element matrix remains the same and as such, the amount of row-locking for each element is constant. In contrast, the p- and d-benchmarks change the size of the element matrix, leading to more row-locking per elememt. Further, the d-benchmark assures contiguous DOF indices on each node, which translates to high vectorization levels. The same is not necessarily true for the p-benchmark, where the new DOF indices are not guaranteed to be contiguous between the element nodes.
\end{remark}

\begin{figure}
	\centering
	\begin{subfigure}{0.48\linewidth}
		\centering
		\includegraphics[width=0.6\linewidth]{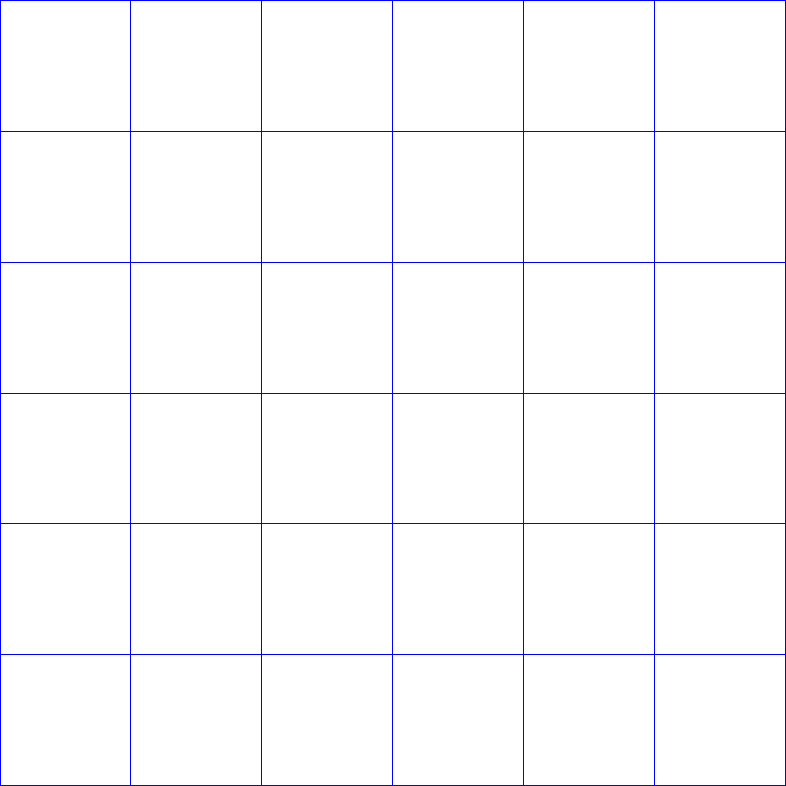}
		\caption{Structured base mesh}
	\end{subfigure}	
	\begin{subfigure}{0.48\linewidth}
		\centering
		\includegraphics[width=0.6\linewidth]{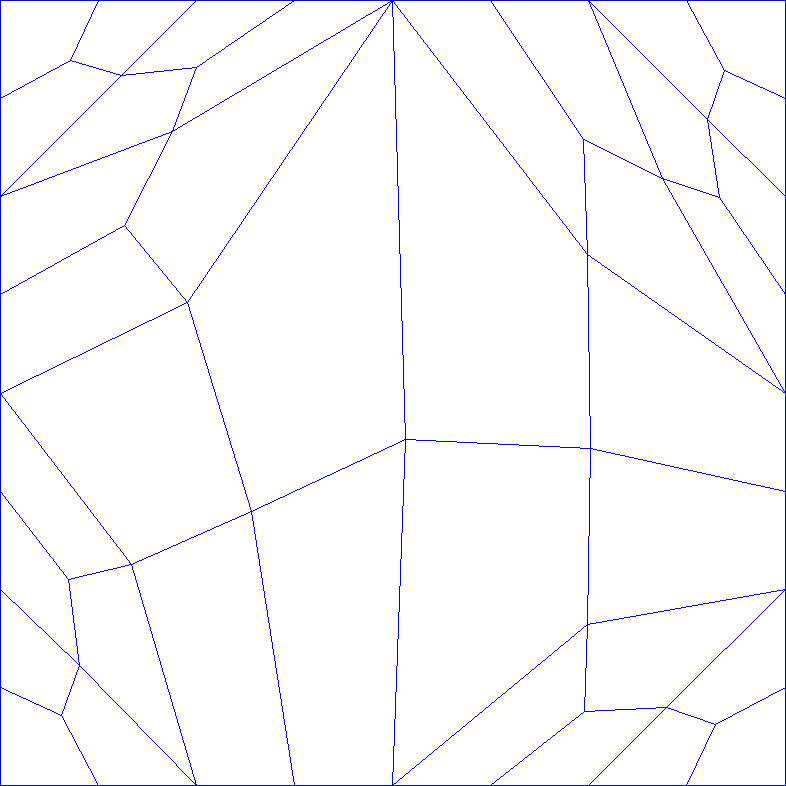}
		\caption{Unstructured base mesh}
	\end{subfigure}	
    \caption{Base benchmark meshes.}
    \label{fig:meshes}
\end{figure}

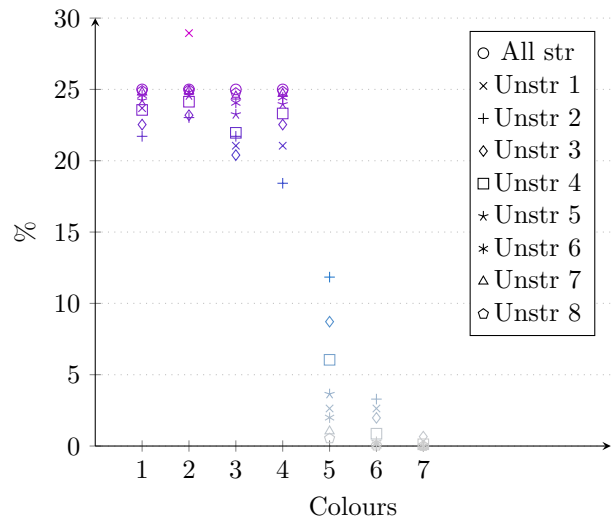
\begin{figure}
	\centering
	\begin{subfigure}{0.48\linewidth}
		\centering
		\begin{tikzpicture}
			\begin{axis}[
				axis lines = left,
				xlabel={Colours},
				ylabel={$\%$},
				xmin=0, xmax=11,
				ymin=0, ymax=30,
				xtick={1,2,3,4,5,6,7},
				ytick={0,5,10,15,20,25,30},
				legend pos=north east,
				ymajorgrids=true,
				grid style=dotted,
				]
				\addplot[
				scatter,
				only marks,
				mark=o,
				]
				coordinates {
					(1,25)
					(2,25)
					(3,25)
					(4,25)
				};
				\addlegendentry{All str}
				\addplot[
				scatter,
				only marks,
				mark=x,
				]
				coordinates {
					(1,9/38*100)
					(2,11/38*100)
					(3,8/38*100)
					(4,8/38*100)
					(5,1/38*100)
					(6,1/38*100)
				};
				\addlegendentry{Unstr 1}
				\addplot[
				scatter,
				only marks,
				mark=+,
				]
				coordinates {
					(1,33/152*100)
					(2,35/152*100)
					(3,33/152*100)
					(4,28/152*100)
					(5,18/152*100)
					(6,5/152*100)
				};
				\addlegendentry{Unstr 2}
				\addplot[
				scatter,
				only marks,
				mark=diamond,
				]
				coordinates {
					(1,137/608*100)
					(2,141/608*100)
					(3,124/608*100)
					(4,137/608*100)
					(5,53/608*100)
					(6,12/608*100)
					(7,4/608*100)
				};
				\addlegendentry{Unstr 3}
				\addplot[
				scatter,
				only marks,
				mark=square,
				]
				coordinates {
					(1,573/2432*100)
					(2,587/2432*100)
					(3,534/2432*100)
					(4,567/2432*100)
					(5,147/2432*100)
					(6,21/2432*100)
					(7,3/2432*100)
				};
				\addlegendentry{Unstr 4}
				\addplot[
				scatter,
				only marks,
				mark=star,
				]
				coordinates {
					(1,2357/9728*100)
					(2,2391/9728*100)
					(3,2262/9728*100)
					(4,2331/9728*100)
					(5,355/9728*100)
					(6,29/9728*100)
					(7,3/9728*100)
				};
				\addlegendentry{Unstr 5}
				\addplot[
				scatter,
				only marks,
				mark=asterisk,
				]
				coordinates {
					(1,9573/38912*100)
					(2,9647/38912*100)
					(3,9362/38912*100)
					(4,9507/38912*100)
					(5,775/38912*100)
					(6,45/38912*100)
					(7,3/38912*100)
				};
				\addlegendentry{Unstr 6}
				\addplot[
				scatter,
				only marks,
				mark=triangle,
				]
				coordinates {
					(1,38597/155648*100)
					(2,38751/155648*100)
					(3,38150/155648*100)
					(4,38451/155648*100)
					(5,1619/155648*100)
					(6,77/155648*100)
					(7,3/155648*100)
				};
				\addlegendentry{Unstr 7}
				\addplot[
				scatter,
				only marks,
				mark=pentagon,
				]
				coordinates {
					(1,155013/622592*100)
					(2,155327/622592*100)
					(3,154090/622592*100)
					(4,154707/622592*100)
					(5,3311/622592*100)
					(6,141/622592*100)
					(7,3/622592*100)
				};
				\addlegendentry{Unstr 8}
			\end{axis}
		\end{tikzpicture}
	\end{subfigure}
	\caption{Distribution of elements per colour for the structured and unstructured meshes, see \cref{fig:meshes}.}
	\label{fig:dis}
\end{figure}

For our experiments we employ several structured and unstructured meshes generated with Gmsh\cite{eigenweb}, starting with a square domain that is progressively refined ,see \cref{fig:meshes}, producing two sets of 8 meshes
with 
\begin{align}
		&\text{structured:} && 36, \, 144, \, 576, \, 2304, \, 9216, \, 36864, \, 147456, \, 589824 \notag \\
		&\text{unstructured:} && 38,\, 152,\, 608, \,2432, \,9728,\, 38912,\, 155648, \,622592 \notag
\end{align}
elements, achieved by splitting each element into $4$. 
Applying the greedy colouring algorithm on the meshes yields the distributions depicted in \cref{fig:dis}. On structured meshes, the algorithm yields the optimal result of $4$ colours with elements distributed evenly. On unstructured meshes, the algorithm finds $6$ to $7$ colours, where $4$ colours contain most of the elements and the remaining $2-3$ colours contain about $10-12 \%$ of the remaining elements.

\begin{remark}
		Every element stiffness matrix is calculated as many times as there are elements in the mesh in each run, as to reflect the equivalent procedure in standard FE-programs.
\end{remark}

For our investigations we employ the following Windows-workstation system
\begin{table}[H]
	\centering
	\begin{tabular}{|l|c|}
		\hline
		\textbf{Processor} & Intel(R) Core(TM) i9-9980HK CPU @ 2.40-5.00 GHz  \\ \hline
		\textbf{Cores}     & 8 Cores                  \\ \hline
		\textbf{AVX up to} & AVX2                                        \\ \hline
		\textbf{RAM}       & 16 GB SODIMM 2667 MHz                          \\ \hline
		\textbf{Compiler}       &  MSVC 16.9.4 - release mode flags: fp:fast, arch:AVX2, O2                          \\ \hline
	\end{tabular}
\end{table}

\begin{figure}
	\begin{subfigure}{0.48\linewidth}
		\centering
		\begin{tikzpicture}
			\pgfplotsset{every tick label/.append style={font=\normalsize},
				every x tick scale label/.append style={yshift=0.3em}}
			\begin{axis}[
				/pgf/number format/1000 sep={},
				axis lines = left,
				xlabel={Global degrees of freedom},
				ylabel={Time $[\mu s]$},
				xmin=0, xmax=7e5,
				ymin=0, ymax= 1.4e5,
				xtick={0, 1e5, 2e5, 3e5, 4e5, 5e5, 6e5, 7e5},
				ytick={0, 2e4, 4e4, 6e4, 8e4, 1e5, 1.2e5, 1.4e5},
				legend pos=north west,
				ymajorgrids=true,
				grid style=dotted,
				]
				\addplot[
				color=blue,
				mark=o,
				]
				coordinates {
					(49,19.39)
					(169,47.09)
					(625,123.35)
					(2401,455.24)
					(9409,1786.51)
					(37249,6862.93)
					(148225,28004.2)
					(591361,117074)
				};
				\addlegendentry{Str Seq $t_\text{avg}$}
				\addplot[
				color=blue,
				style=dashed,
				mark=x,
				]
				coordinates {
					(49,13)
					(169,35)
					(625,115)
					(2401,448)
					(9409,1756)
					(37249,6504)
					(148225,26526)
					(591361,112029)
				};
				\addlegendentry{Str Seq $t_\text{min}$}
				\addplot[
				color=teal,
				mark=square,
				]
				coordinates {
					(55,19.75)
					(185,55.91)
					(673,145.1)
					(2561,545.87)
					(9985,2055.54)
					(39425,7764.97)
					(156673,31682.6)
					(624641,122767)
				};
				\addlegendentry{Unstr Seq $t_\text{avg}$}
				\addplot[
				color=teal,
				style=dashed,
				mark=star,
				]
				coordinates {
					(55,14)
					(185,40)
					(673,135)
					(2561,538)
					(9985,1966)
					(39425,7450)
					(156673,30184)
					(624641,119676)
				};
				\addlegendentry{Unstr Seq $t_\text{min}$}
			\end{axis}
		\end{tikzpicture}
		\caption{h-benchmark}
		\label{fig:seqh}
	\end{subfigure}
	\begin{subfigure}{0.48\linewidth}
		\centering
		\begin{tikzpicture}
			\pgfplotsset{every tick label/.append style={font=\normalsize},
				every x tick scale label/.append style={yshift=0.3em}}
			\begin{axis}[
				/pgf/number format/1000 sep={},
				axis lines = left,
				xlabel={Global degrees of freedom},
				ylabel={Time $[\mu s]$},
				xmin=0, xmax=160000,
				ymin=0, ymax=240000,
				xtick={0,40000,80000,120000, 160000},
				ytick={0, 40e3, 80e3,120e3,160e3,200e3, 240e3},
				legend pos=north west,
				ymajorgrids=true,
				grid style=dotted,
				]
				\addplot[
				color=blue,
				mark=o,
				]
				coordinates {
					(2401,469.78)
					(9409,2784.85)
					(21025,10215.8)
					(37249,19335)
					(58081,38796.3)
					(83521,78211.9)
					(113569,130165)
					(148225,225772)
				};
				\addlegendentry{Str Seq $t_\text{avg}$}
				\addplot[
				color=blue,
				style=dashed,
				mark=x,
				]
				coordinates {
					(2401,459)
					(9409,2766)
					(21025,9841)
					(37249,18603)
					(58081,37169)
					(83521,75561)
					(113569,127098)
					(148225,220950)
				};
				\addlegendentry{Str Seq $t_\text{min}$}
				\addplot[
				color=teal,
				mark=square,
				]
				coordinates {
					(2561,527.54)
					(9985,3128.2)
					(22273,11118.8)
					(39425,32897.4)
					(61441,41407.9)
					(88321,82590.8)
					(120065,137704)
					(156673,239250)
				};
				\addlegendentry{Unstr Seq $t_\text{avg}$}
				\addplot[
				color=teal,
				style=dashed,
				mark=star,
				]
				coordinates {
					(2561,520)
					(9985,3109)
					(22273,10584)
					(39425,31596)
					(61441,39370)
					(88321,80091)
					(120065,135376)
					(156673,234772)
				};
				\addlegendentry{Unstr Seq $t_\text{min}$}
			\end{axis}
		\end{tikzpicture}
		\caption{p-benchmark}
		\label{fig:seqp}
	\end{subfigure}
	\caption{Assembly times of the sequential method for a scalar problem $d = 1$.}
	\label{fig:seq}
\end{figure}
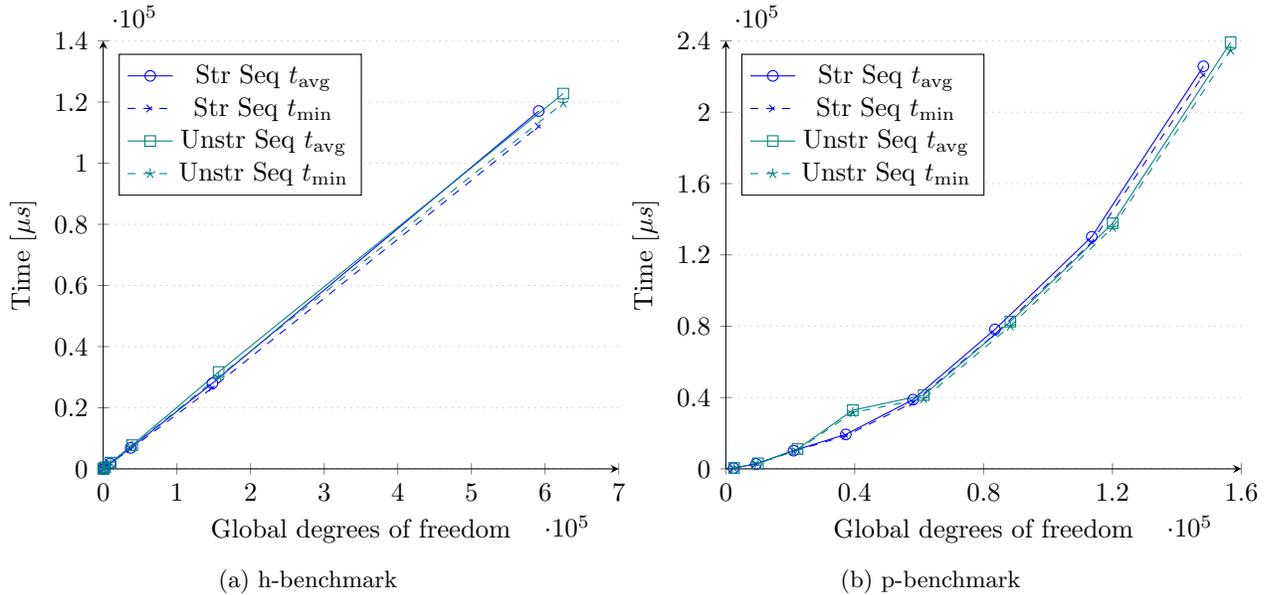

The performance of the sequential assembly method for the case $d = 1$ is given in \cref{fig:seq} as a frame of reference for computation times, where the abbreviations Str and Unstr refer to structured and unstructured meshes, respectively.
The computation time of the h-benchmark grows linearly, whereas the slope of the computation times of the p-benchmark for the fourth structured and unstructured meshes changes for each point. This can be traced back to the ever growing size of the element matrix and consequently, the higher amount of time required to initialize and assemble each element matrix.

The number of none zero entries for the highest amount of DOF for each benchmark reads

\begin{table}[H]
	\centering
	\begin{tabular}{|l|c|c|}
		\hline
		h-benchmark & Str: 5313025 & Unstr: 5609473 \\ \hline
		p-benchmark  & Str: 14753281  & Unstr: 15575041  \\ \hline
	\end{tabular}
\end{table}

The performance of the sequenced assembly method for the case $d = 4$ is depicted in \cref{fig:seqd4} and behaves similarly to the case $d=1$. The number of none zero entries for the highest amount of DOF for each of the benchmarks for the case $d = 4$ is given by

\begin{table}[H]
	\centering
	\begin{tabular}{|l|c|c|}
		\hline
		h-benchmark & Str: 5326864 & Unstr: 5627920 \\ \hline
		p-benchmark  & Str: 21270544  & Unstr: 22462480  \\ \hline
	\end{tabular}
\end{table}

\begin{figure}
	\centering
	\begin{subfigure}{0.48\linewidth}
		\centering
		\begin{tikzpicture}
			\pgfplotsset{every tick label/.append style={font=\normalsize},
				every x tick scale label/.append style={yshift=0.3em}}
			\begin{axis}[
				/pgf/number format/1000 sep={},
				axis lines = left,
				xlabel={Global degrees of freedom},
				ylabel={Time $[\mu s]$},
				xmin=0, xmax=320e3,
				ymin=0, ymax=600e3,
				xtick={0,80e3,160e3,240e3, 320e3},
				ytick={0,100e3,200e3,300e3, 400e3, 500e3, 600e3},
				legend pos=north west,
				ymajorgrids=true,
				grid style=dotted,
				]
				\addplot[
				color=blue,
				mark=o,
				]
				coordinates {
					(37249,7076.59)
					(74498,31396.9)
					(111747,69706.5)
					(148996,128901)
					(186245,209184)
					(223494,305040)
					(260743,372297)
					(297992,557512)
				};
				\addlegendentry{Str Seq $t_\text{avg}$}
				\addplot[
				color=blue,
				style=dashed,
				mark=x,
				]
				coordinates {
					(37249,6610)
					(74498,30146)
					(111747,66605)
					(148996,125402)
					(186245,204200)
					(223494,297099)
					(260743,363792)
					(297992,539043)
				};
				\addlegendentry{Str Seq $t_\text{min}$}
				\addplot[
				color=teal,
				mark=square,
				]
				coordinates {
					(39425,7523.68)
					(78850,34071.8)
					(118275,74094)
					(157700,136654)
					(197125,220003)
					(236550,319044)
					(275975,392885)
					(315400,580393)
				};
				\addlegendentry{Unstr Seq $t_\text{avg}$}
				\addplot[
				color=teal,
				style=dashed,
				mark=star,
				]
				coordinates {
					(39425,7312)
					(78850,32367)
					(118275,71662)
					(157700,133850)
					(197125,215886)
					(236550,311135)
					(275975,384412)
					(315400,566696)
				};
				\addlegendentry{Unstr Seq $t_\text{min}$}
			\end{axis}
		\end{tikzpicture}
	\end{subfigure}
	\caption{Assembly times of the d-benchmark for the sequential method.}
	\label{fig:seqd}
\end{figure}
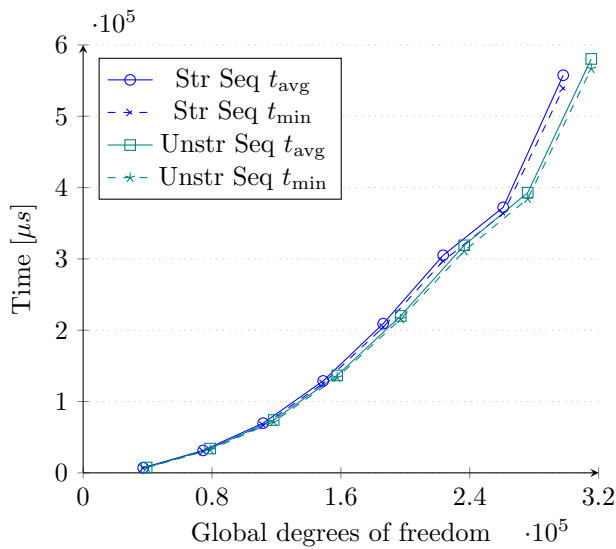

Further, we give the assembly times of the d-benchmark on the sixth structured and unstructured meshes, see \cref{fig:seqd}. where the highest numbers of non-zero entries read

\begin{table}[H]
	\centering
	\begin{tabular}{|l|c|c|}
		\hline
		d-benchmark & Str: 21307456 & Unstr: 22511680 \\ \hline
	\end{tabular}
\end{table}

\begin{remark}
		The size of the values array is given by $s = (\nicefrac{8}{1000000}) \, n_\text{entires} $ in megabytes.
		Consequently, for value arrays with more than $2000000$ entries the array cannot be contained in the L3 cache, whose limit is 16 megabytes. However, this is a simplification. In reality, even smaller problems might need to access the RAM, as other data structures such the rows and columns arrays or the elements array are also being loaded into the cache.
\end{remark}

\begin{figure}
	\begin{subfigure}{0.48\linewidth}
		\centering
		\begin{tikzpicture}
			\pgfplotsset{every tick label/.append style={font=\normalsize},
				every x tick scale label/.append style={yshift=0.3em}}
			\begin{axis}[
				/pgf/number format/1000 sep={},
				axis lines = left,
				xlabel={Global degrees of freedom},
				ylabel={Time $[\mu s]$},
				xmin=0, xmax=1.6e5,
				ymin=0, ymax= 1.4e5,
				xtick={0,0.4e5, 0.8e5, 1.2e5, 1.6e5},
				ytick={0, 2e4, 4e4, 6e4, 8e4, 1e5, 1.2e5, 1.4e5},
				legend pos=north west,
				ymajorgrids=true,
				grid style=dotted,
				]
				\addplot[
				color=blue,
				mark=o,
				]
				coordinates {
					(196,133.53)
					(676,490.05)
					(2500,1986.47)
					(9604,7727.24)
					(37636,31054.8)
					(148996,128376)
				};
				\addlegendentry{Str Seq $t_\text{avg}$}
				\addplot[
				color=blue,
				style=dashed,
				mark=x,
				]
				coordinates {
					(196,124)
					(676,488)
					(2500,1869)
					(9604,7267)
					(37636,30011)
					(148996,124853)
				};
				\addlegendentry{Str Seq $t_\text{min}$}
				\addplot[
				color=teal,
				mark=square,
				]
				coordinates {
					(220,182.62)
					(740,536.79)
					(2692,2164.64)
					(10244,8316.04)
					(39940,33425.2)
					(157700,136257)
				};
				\addlegendentry{Unstr Seq $t_\text{avg}$}
				\addplot[
				color=teal,
				style=dashed,
				mark=star,
				]
				coordinates {
					(220,179)
					(740,530)
					(2692,2107)
					(10244,8037)
					(39940,33068)
					(157700,133746)
				};
				\addlegendentry{Unstr Seq $t_\text{min}$}
			\end{axis}
		\end{tikzpicture}
		\caption{h-benchmark}
		\label{fig:seqd4h}
	\end{subfigure}
	\begin{subfigure}{0.48\linewidth}
		\centering
		\begin{tikzpicture}
			\pgfplotsset{every tick label/.append style={font=\normalsize},
				every x tick scale label/.append style={yshift=0.3em}}
			\begin{axis}[
				/pgf/number format/1000 sep={},
				axis lines = left,
				xlabel={Global degrees of freedom},
				ylabel={Time $[\mu s]$},
				xmin=0, xmax=160000,
				ymin=0, ymax=420000,
				xtick={0,40000,80000,120000, 160000},
				ytick={0, 105e3, 210e3, 315e3,420e3},
				legend pos=north west,
				ymajorgrids=true,
				grid style=dotted,
				]
				\addplot[
				color=blue,
				mark=o,
				]
				coordinates {
					(9604,7463.98)
					(37636,43359)
					(84100,144142)
					(148996,389784)
				};
				\addlegendentry{Str Seq $t_\text{avg}$}
				\addplot[
				color=blue,
				style=dashed,
				mark=x,
				]
				coordinates {
					(9604,7210)
					(37636,41615)
					(84100,140552)
					(148996,384223)
				};
				\addlegendentry{Str Seq $t_\text{min}$}
				\addplot[
				color=teal,
				mark=square,
				]
				coordinates {
					(10244,8262.94)
					(39940,46409)
					(89092,154510)
					(157700,411737)
				};
				\addlegendentry{Unstr Seq $t_\text{avg}$}
				\addplot[
				color=teal,
				style=dashed,
				mark=star,
				]
				coordinates {
					(10244,7725)
					(39940,44717)
					(89092,150730)
					(157700,404403)
				};
				\addlegendentry{Unstr Seq $t_\text{min}$}
			\end{axis}
		\end{tikzpicture}
		\caption{p-benchmark}
		\label{fig:seqd4p}
	\end{subfigure}
	\caption{Assembly times of the sequential method for a vector problem $d = 4$.}
	\label{fig:seqd4}
\end{figure}
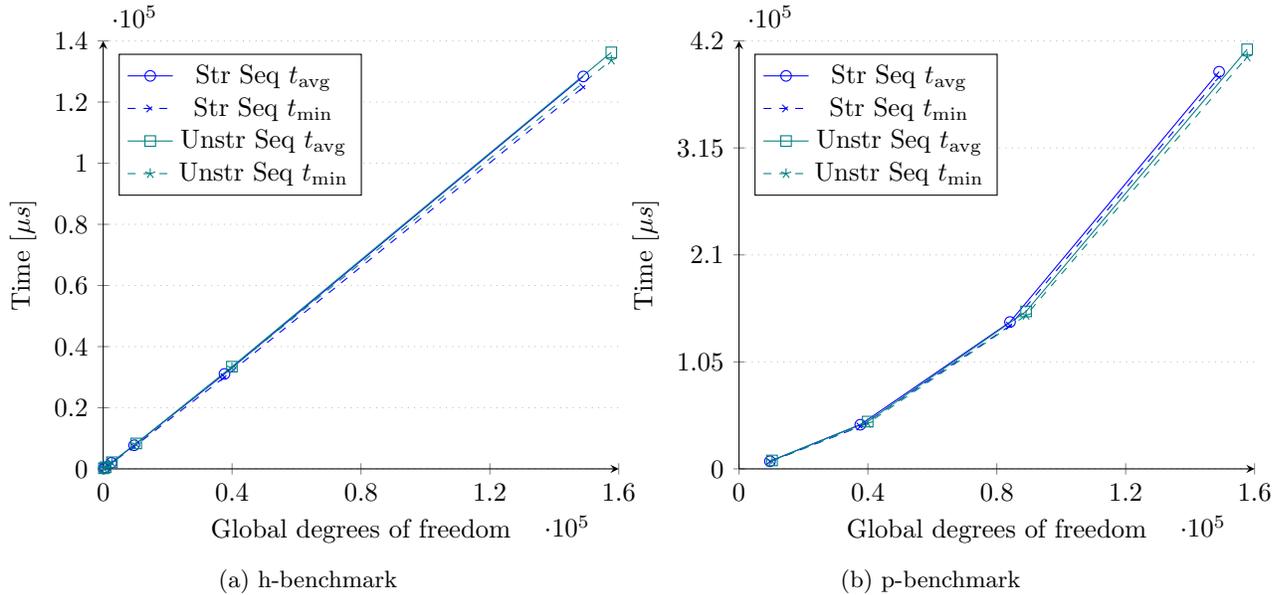

Finally, we compare the storage requirements of the CSR and CRAC formats. To that end, we define the factor
\begin{align}
	\gamma = \dfrac{n^\text{CRAC}_\text{ca}}{n^\text{CSR}_\text{ci}} \, ,
\end{align}
where $n^\text{CSR}_\text{ci}$ represents the length of the column indices array of the CSR format and $n^\text{CRAC}_\text{ca}$ the length of the column alignments array of the CRAC format. For $\gamma < 1$ the CRAC format represents the lighter storage scheme and for $\gamma > 1$ the heavier, since both formats have the same length for the row pointers and value arrays.

For all benchmark figures, the methods are summarized as follows:

\begin{table}[H]
	\centering
	\begin{tabular}{|c|c|c|c|c|c|}
		\hline
		\textbf{Method:}  & \textbf{Abbreviation:} & \textbf{Parallel:} & \textbf{Vectorized:} & \textbf{for CSR:} & \textbf{for CRAC:} \\ \hline &&&&& \\[-2ex] 
		 Atomic \lstinline|fetch_add| & $Atc^\text{CSR}$ & $\surd$ & & $\surd$& \\ \hline &&&&& \\[-2ex] 
		 Using \lstinline|spin_int| & $Sp^\text{CSR}$ & $\surd$ & & $\surd$& \\ \hline &&&&& \\[-2ex] 
		 Using colouring & $Col^\text{CSR}_\text{vec}$ & $\surd$ & $\surd$ & $\surd$& \\ \hline &&&&& \\[-2ex] 
		 Vectorized with \lstinline|spin_int| & $Sp_\text{vec}^\text{CSR}, \quad Sp_\text{vec}^\text{CRAC}$  & $\surd$ & $\surd$ & $\surd$& $\surd$ \\ \hline
	\end{tabular}
\end{table}

\subsection{Speed factors}
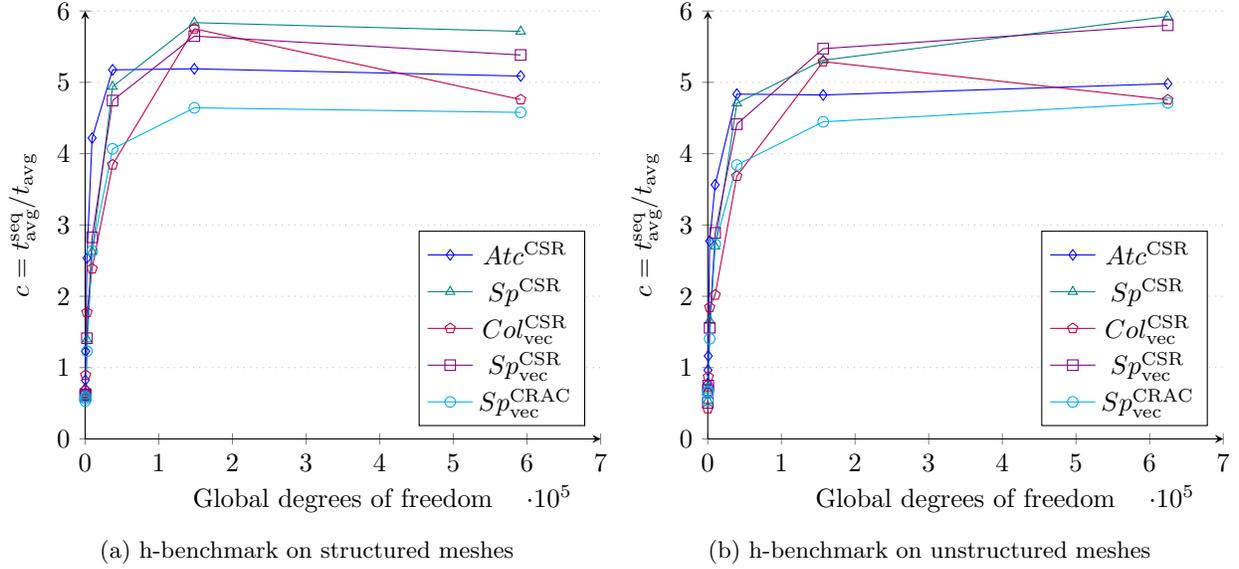
\begin{figure}
	\begin{subfigure}{0.48\linewidth}
		\centering
		\begin{tikzpicture}
			\pgfplotsset{every tick label/.append style={font=\normalsize},
				every x tick scale label/.append style={yshift=0.3em}}
			\begin{axis}[
				/pgf/number format/1000 sep={},
				axis lines = left,
				xlabel={Global degrees of freedom},
				ylabel={$c = t_\text{avg}^\text{seq} / t_\text{avg}$},
				xmin=0, xmax=7e5,
				ymin=0, ymax= 6,
				xtick={0, 1e5, 2e5, 3e5, 4e5, 5e5, 6e5, 7e5},
				ytick={0,1,2,3,4,5,6},
				legend pos=south east,
				ymajorgrids=true,
				grid style=dotted,
				]
				\addplot[
				color=blue,
				mark=diamond,
				]
				coordinates {
					(49,19.39/27.49)
					(169,47.09/57.11)
					(625,123.35/100.42)
					(2401,455.24/179.51)
					(9409,1786.51/423.55)
					(37249,6862.93/1326.07)
					(148225,28004.2/5396.39)
					(591361,117074/23009.7)
				};
				\addlegendentry{$Atc^\text{CSR}$}
				\addplot[
				color=teal,
				mark=triangle,
				]
				coordinates {
					(49,19.39/32.01)
					(169,47.09/79.74)
					(625,123.35/199.11)
					(2401,455.24/325.62)
					(9409,1786.51/680.49)
					(37249,6862.93/1389.32)
					(148225,28004.2/4798)
					(591361,117074/20489.3)
				};
				\addlegendentry{$Sp^\text{CSR}$}
				\addplot[
				color=purple,
				mark=pentagon,
				]
				coordinates {
					(49,19.39/28.58)
					(169,47.09/73.04)
					(625,123.35/137.98)
					(2401,455.24/256.52)
					(9409,1786.51/748.43)
					(37249,6862.93/1785.41)
					(148225,28004.2/4868.91)
					(591361,117074/24604.8)
				};
				\addlegendentry{$Col^\text{CSR}_\text{vec}$}
				\addplot[
				color=violet,
				mark=square,
				]
				coordinates {
					(49,19.39/31.4)
					(169,47.09/78.59)
					(625,123.35/198.68)
					(2401,455.24/322.46)
					(9409,1786.51/631.9)
					(37249,6862.93/1446.54)
					(148225,28004.2/4957.88)
					(591361,117074/21746.6)
				};
				\addlegendentry{$Sp_\text{vec}^\text{CSR}$}
				\addplot[
				color=cyan,
				mark=o,
				]
				coordinates {
					(49,19.39/36.53)
					(169,47.09/83.41)
					(625,123.35/207.93)
					(2401,455.24/370.2)
					(9409,1786.51/677.98)
					(37249,6862.93/1687.66)
					(148225,28004.2/6032.21)
					(591361,117074/25572.2)
				};
				\addlegendentry{$Sp_\text{vec}^\text{CRAC}$}
			\end{axis}
		\end{tikzpicture}
		\caption{h-benchmark on structured meshes}
	\end{subfigure}
	\begin{subfigure}{0.48\linewidth}
		\centering
		\begin{tikzpicture}
			\pgfplotsset{every tick label/.append style={font=\normalsize},
				every x tick scale label/.append style={yshift=0.3em}}
			\begin{axis}[
				/pgf/number format/1000 sep={},
				axis lines = left,
				xlabel={Global degrees of freedom},
				ylabel={$c = t_\text{avg}^\text{seq} / t_\text{avg}$},
				xmin=0, xmax=7e5,
				ymin=0, ymax= 6,
				xtick={0, 1e5, 2e5, 3e5, 4e5, 5e5, 6e5, 7e5},
				ytick={0,1,2,3,4,5,6},
				legend pos=south east,
				ymajorgrids=true,
				grid style=dotted,
				]
				\addplot[
				color=blue,
				mark=diamond,
				]
				coordinates {
					(55,19.75/25.39)
					(185,55.91/58.05)
					(673,145.1/124.73)
					(2561,545.87/196.71)
					(9985,2055.54/576.91)
					(39425,7764.97/1606.15)
					(156673,31682.6/6570.01)
					(624641,122767/24648.8)
				};
				\addlegendentry{$Atc^\text{CSR}$}
				\addplot[
				color=teal,
				mark=triangle,
				]
				coordinates {
					(55,19.75/36.5)
					(185,55.91/81.82)
					(673,145.1/195.23)
					(2561,545.87/329.03)
					(9985,2055.54/760.9)
					(39425,7764.97/1649.48)
					(156673,31682.6/5967.01)
					(624641,122767/20722.6)
				};
				\addlegendentry{$Sp^\text{CSR}$}
				\addplot[
				color=purple,
				mark=pentagon,
				]
				coordinates {
					(55,19.75/46.44)
					(185,55.91/86.14)
					(673,145.1/166.83)
					(2561,545.87/295.71)
					(9985,2055.54/1018.24)
					(39425,7764.97/2107.64)
					(156673,31682.6/5989.12)
					(624641,122767/25816.5)
				};
				\addlegendentry{$Col^\text{CSR}_\text{vec}$}
				\addplot[
				color=violet,
				mark=square,
				]
				coordinates {
					(55,19.75/39.02)
					(185,55.91/81.71)
					(673,145.1/193.92)
					(2561,545.87/350.17)
					(9985,2055.54/711.47)
					(39425,7764.97/1759.35)
					(156673,31682.6/5789.95)
					(624641,122767/21167.9)
				};
				\addlegendentry{$Sp_\text{vec}^\text{CSR}$}
				\addplot[
				color=cyan,
				mark=o,
				]
				coordinates {
					(55,19.75/36.93)
					(185,55.91/85.66)
					(673,145.1/216.89)
					(2561,545.87/388.81)
					(9985,2055.54/752.61)
					(39425,7764.97/2021.81)
					(156673,31682.6/7125.84)
					(624641,122767/26052.1)
				};
				\addlegendentry{$Sp_\text{vec}^\text{CRAC}$}
			\end{axis}
		\end{tikzpicture}
		\caption{h-benchmark on unstructured meshes}
	\end{subfigure}
	\caption{Improvement factors of assembly times for the h-benchmark for the case $d = 1$.}
	\label{fig:par}
\end{figure}

We start this discussion by considering the case $d = 1$. The speed factors for the h-benchmark are given in \cref{fig:par}. We notice the slopes of the improvement curves flatten with ever larger meshes. Generally, for all parallel methods we find an improvement factor between $4.5$ and $6$, representing the limit for improvement for the parallel methods for this type of problems on our system. For both mesh types on the finest mesh, the fastest method was the non-vectorized \lstinline|spin_int| method, yielding the factors $c = 5.71$ and $c = 5.92$ for the structured and unstructured meshes, respectively. Since both the colouring and atomic methods are lock-free, it is clear that the locking algorithm \textbf{does not} dominate the assembly time in the \lstinline|spin_int| methods. Further, it suggests lock-free methods are not necessarily optimal for parallel assembly as, although the atomic method allows parallel access to all elements of the sparse matrix, it falls short in terms of performance in comparison with the \lstinline|spin_int| method. In both cases the CRAC format performed slightly slower than the CSR format. This is likely due to the low amount of contiguous columns present in the $d = 1$ case, making the CRAC format the heavier storage scheme.

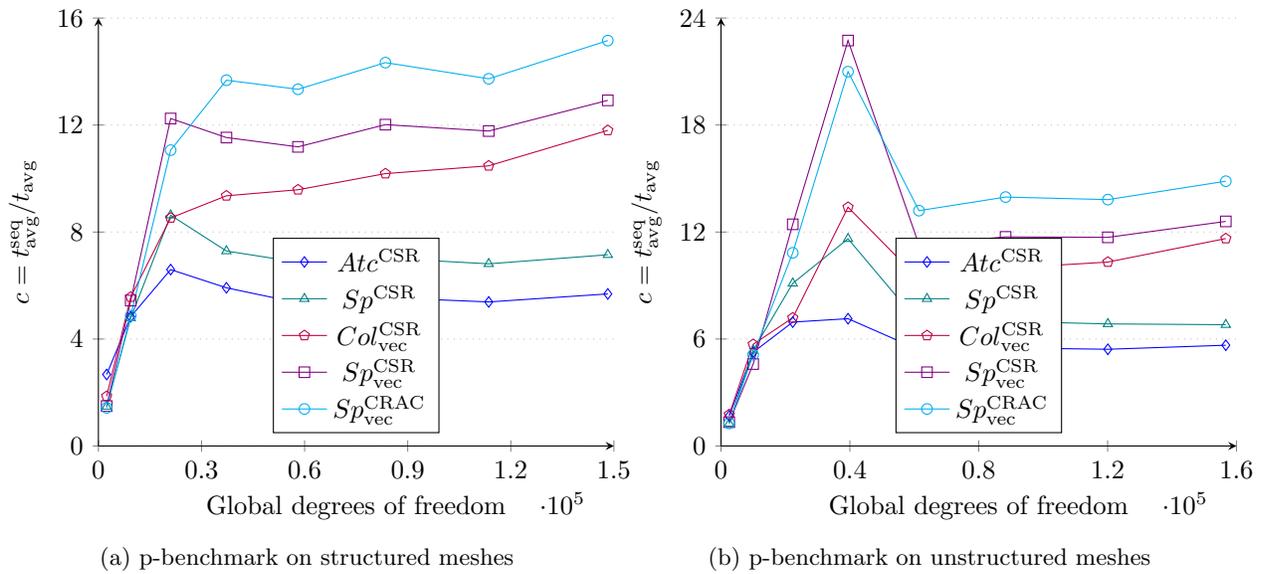
\begin{figure}
	\begin{subfigure}{0.48\linewidth}
		\centering
		\begin{tikzpicture}
			\pgfplotsset{every tick label/.append style={font=\normalsize},
				every x tick scale label/.append style={yshift=0.3em}}
			\begin{axis}[
				/pgf/number format/1000 sep={},
				axis lines = left,
				xlabel={Global degrees of freedom},
				ylabel={$c = t_\text{avg}^\text{seq} / t_\text{avg}$},
				xmin=0, xmax=150e3,
				ymin=0, ymax= 16,
				xtick={0, 30e3, 60e3, 90e3, 120e3, 150e3},
				ytick={0,4,8,12,16},
				legend style={at={(0.5,0.03)},anchor=south},
				ymajorgrids=true,
				grid style=dotted,
				]
				\addplot[
				color=blue,
				mark=diamond,
				]
				coordinates {
					(2401,469.78/175.67)
					(9409,2784.85/575.8)
					(21025,10215.8/1548.83)
					(37249,19335/3268.97)
					(58081,38796.3/7247.2)
					(83521,78211.9/14095.1)
					(113569,130165/24166.2)
					(148225,225772/39703.8)
				};
				\addlegendentry{$Atc^\text{CSR}$}
				\addplot[
				color=teal,
				mark=triangle,
				]
				coordinates {
					(2401,469.78/318)
					(9409,2784.85/581.79)
					(21025,10215.8/1183.8)
					(37249,19335/2652.19)
					(58081,38796.3/5645.7)
					(83521,78211.9/11134)
					(113569,130165/19108.4)
					(148225,225772/31573.5)
				};
				\addlegendentry{$Sp^\text{CSR}$}
				\addplot[
				color=purple,
				mark=pentagon,
				]
				coordinates {
					(2401,469.78/252.58)
					(9409,2784.85/500.08)
					(21025,10215.8/1197.86)
					(37249,19335/2066.73)
					(58081,38796.3/4049.3)
					(83521,78211.9/7679.77)
					(113569,130165/12424.9)
					(148225,225772/19133.5)
				};
				\addlegendentry{$Col^\text{CSR}_\text{vec}$}
				\addplot[
				color=violet,
				mark=square,
				]
				coordinates {
					(2401,469.78/314.38)
					(9409,2784.85/512.06)
					(21025,10215.8/834.41)
					(37249,19335/1676.54)
					(58081,38796.3/3468.8)
					(83521,78211.9/6509.14)
					(113569,130165/11054)
					(148225,225772/17469.8)
				};
				\addlegendentry{$Sp_\text{vec}^\text{CSR}$}
				\addplot[
				color=cyan,
				mark=o,
				]
				coordinates {
					(2401,469.78/329.19)
					(9409,2784.85/571.47)
					(21025,10215.8/923.61)
					(37249,19335/1413.67)
					(58081,38796.3/2909.04)
					(83521,78211.9/5456.2)
					(113569,130165/9480.65)
					(148225,225772/14894.2)
				};
				\addlegendentry{$Sp_\text{vec}^\text{CRAC}$}
			\end{axis}
		\end{tikzpicture}
		\caption{p-benchmark on structured meshes}
	\end{subfigure}
	\begin{subfigure}{0.48\linewidth}
		\centering
		\begin{tikzpicture}
			\pgfplotsset{every tick label/.append style={font=\normalsize},
				every x tick scale label/.append style={yshift=0.3em}}
			\begin{axis}[
				/pgf/number format/1000 sep={},
				axis lines = left,
				xlabel={Global degrees of freedom},
				ylabel={$c = t_\text{avg}^\text{seq} / t_\text{avg}$},
				xmin=0, xmax=160e3,
				ymin=0, ymax= 24,
				xtick={0, 40e3, 80e3, 120e3, 160e3},
				ytick={0,6,12,18,24},
				legend style={at={(0.5,0.03)},anchor=south},
				ymajorgrids=true,
				grid style=dotted,
				]
				\addplot[
				color=blue,
				mark=diamond,
				]
				coordinates {
					(2561,527.54/312.59)
					(9985,3128.2/594.89)
					(22273,11118.8/1599.94)
					(39425,32897.4/4607.65)
					(61441,41407.9/7674.81)
					(88321,82590.8/15011.3)
					(120065,137704/25383)
					(156673,239250/42331.6)
				};
				\addlegendentry{$Atc^\text{CSR}$}
				\addplot[
				color=teal,
				mark=triangle,
				]
				coordinates {
					(2561,527.54/417.66)
					(9985,3128.2/582.56)
					(22273,11118.8/1218.84)
					(39425,32897.4/2830.66)
					(61441,41407.9/6000.73)
					(88321,82590.8/11815.5)
					(120065,137704/20110.7)
					(156673,239250/35167.5)
				};
				\addlegendentry{$Sp^\text{CSR}$}
				\addplot[
				color=purple,
				mark=pentagon,
				]
				coordinates {
					(2561,527.54/298.92)
					(9985,3128.2/548.66)
					(22273,11118.8/1546.54)
					(39425,32897.4/2456.56)
					(61441,41407.9/4424.56)
					(88321,82590.8/8283.6)
					(120065,137704/13348.4)
					(156673,239250/20563.6)
				};
				\addlegendentry{$Col^\text{CSR}_\text{vec}$}
				\addplot[
				color=violet,
				mark=square,
				]
				coordinates {
					(2561,527.54/394.43)
					(9985,3128.2/680.44)
					(22273,11118.8/894.28)
					(39425,32897.4/1446.68)
					(61441,41407.9/3668.85)
					(88321,82590.8/7041.57)
					(120065,137704/11767.2)
					(156673,239250/18994.6)
				};
				\addlegendentry{$Sp_\text{vec}^\text{CSR}$}
				\addplot[
				color=cyan,
				mark=o,
				]
				coordinates {
					(2561,527.54/416.49)
					(9985,3128.2/611.23)
					(22273,11118.8/1026.34)
					(39425,32897.4/1566.93)
					(61441,41407.9/3137.13)
					(88321,82590.8/5915.95)
					(120065,137704/9965.78)
					(156673,239250/16111.9)
				};
				\addlegendentry{$Sp_\text{vec}^\text{CRAC}$}
			\end{axis}
		\end{tikzpicture}
		\caption{p-benchmark on unstructured meshes}
	\end{subfigure}
	\caption{Improvement factors of assembly times for the p-benchmark for the case $d = 1$ using the fourth structured and unstructured meshes.}
	\label{fig:par_p}
\end{figure}

The performance of the p-benchmark using the fourth structured and unstructured meshes is depicted in \cref{fig:par_p}. We notice a sudden change in the speed improvement factor between 39425 and 61441 degrees of freedom on unstructured meshes. At 39425 DOF, the size of the values array is about $11.2$ megabytes. At 61441 DOF, the size changes to about $23.2$ megabytes. In other words, up to and including 39425 DOF, access to the sparse matrix was done purely on the cache. From 61441 DOF and on, the matrix no longer fits in the L3 cache($16$ megabytes) and RAM access is required, being much slower. Interestingly, the phenomenon is less apparent for structured meshes. This suggests the structured nature of the meshes facilitates memory fetching, which in turn, allows to slightly diminish the negative effects of the memory wall.
The results of the p-benchmark reinforce the conclusions of the h-benchmark. Namely, using lock-free methods does not guarantee higher performance. Additionally, the CRAC format now performs better than the CSR format for $p > 3$. This is a direct result of the naturally occurring, higher level of contiguity of columns in the sparse matrix. Another consequence of the p-adaption is an increased level of contiguity of column index sequences on each element and the resulting higher performance of the vectorized methods in comparison with the non-vectorized methods.  

\begin{remark}
		Note that the higher level of contiguity of the column alignments array in the CRAC format is not hard-coded in any way. In fact, it results naturally from the higher order elements. Further, the elements produce partially contiguous DOF indices. 
\end{remark}

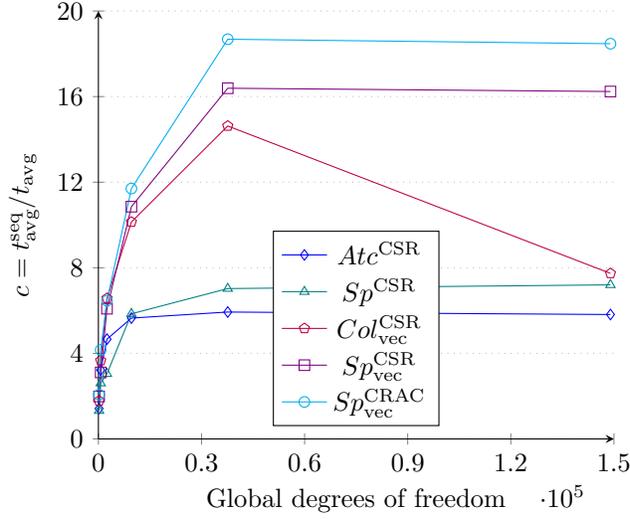
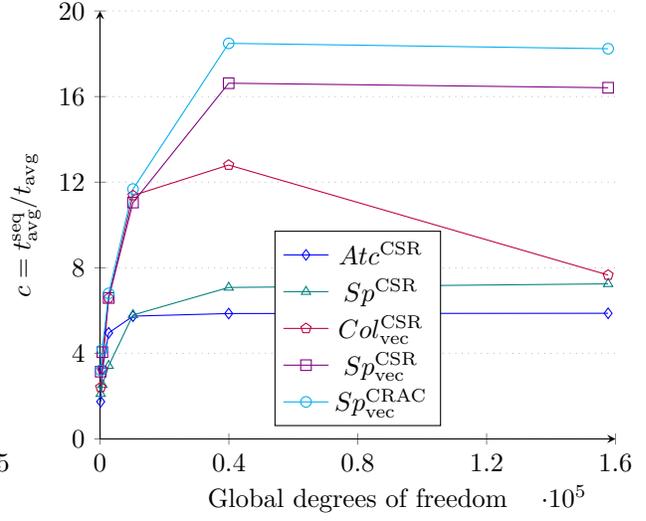
\begin{figure}
	\begin{subfigure}{0.48\linewidth}
		\centering
		\begin{tikzpicture}
			\pgfplotsset{every tick label/.append style={font=\normalsize},
				every x tick scale label/.append style={yshift=0.3em}}
			\begin{axis}[
				/pgf/number format/1000 sep={},
				axis lines = left,
				xlabel={Global degrees of freedom},
				ylabel={$c = t_\text{avg}^\text{seq} / t_\text{avg}$},
				xmin=0, xmax=150e3,
				ymin=0, ymax= 20,
				xtick={0, 30e3, 60e3, 90e3, 120e3, 150e3},
				ytick={0,4,8,12,16, 20},
				legend style={at={(0.5,0.03)},anchor=south},
				ymajorgrids=true,
				grid style=dotted,
				]
				\addplot[
				color=blue,
				mark=diamond,
				]
				coordinates {
					(196,133.53/94.11)
					(676,490.05/152.06)
					(2500,1986.47/425.95)
					(9604,7727.24/1368.03)
					(37636,31054.8/5236.18)
					(148996,128376/22075.6)
				};
				\addlegendentry{$Atc^\text{CSR}$}
				\addplot[
				color=teal,
				mark=triangle,
				]
				coordinates {
					(196,133.53/100.97)
					(676,490.05/187.96)
					(2500,1986.47/651.44)
					(9604,7727.24/1321.19)
					(37636,31054.8/4419.02)
					(148996,128376/17826.8)
				};
				\addlegendentry{$Sp^\text{CSR}$}
				\addplot[
				color=purple,
				mark=pentagon,
				]
				coordinates {
					(196,133.53/75.02)
					(676,490.05/133.76)
					(2500,1986.47/302.18)
					(9604,7727.24/762.21)
					(37636,31054.8/2123.21)
					(148996,128376/16592.5)
				};
				\addlegendentry{$Col^\text{CSR}_\text{vec}$}
				\addplot[
				color=violet,
				mark=square,
				]
				coordinates {
					(196,133.53/67.23)
					(676,490.05/157.63)
					(2500,1986.47/326.37)
					(9604,7727.24/711.91)
					(37636,31054.8/1894.77)
					(148996,128376/7906.17)
				};
				\addlegendentry{$Sp_\text{vec}^\text{CSR}$}
				\addplot[
				color=cyan,
				mark=o,
				]
				coordinates {
					(196,133.53/64.51)
					(676,490.05/117.77)
					(2500,1986.47/307.44)
					(9604,7727.24/660.28)
					(37636,31054.8/1662.42)
					(148996,128376/6949.88)
				};
				\addlegendentry{$Sp_\text{vec}^\text{CRAC}$}
			\end{axis}
		\end{tikzpicture}
		\caption{h-benchmark on structured meshes}
	\end{subfigure}
	\begin{subfigure}{0.48\linewidth}
		\centering
		\begin{tikzpicture}
			\pgfplotsset{every tick label/.append style={font=\normalsize},
				every x tick scale label/.append style={yshift=0.3em}}
			\begin{axis}[
				/pgf/number format/1000 sep={},
				axis lines = left,
				xlabel={Global degrees of freedom},
				ylabel={$c = t_\text{avg}^\text{seq} / t_\text{avg}$},
				xmin=0, xmax=160e3,
				ymin=0, ymax= 20,
				xtick={0, 40e3, 80e3, 120e3, 160e3},
				ytick={0,4,8,12,16,20},
				legend style={at={(0.5,0.03)},anchor=south},
				ymajorgrids=true,
				grid style=dotted,
				]
				\addplot[
				color=blue,
				mark=diamond,
				]
				coordinates {
					(220,182.62/104.84)
					(740,536.79/165)
					(2692,2164.64/436.71)
					(10244,8316.04/1448.97)
					(39940,33425.2/5703.26)
					(157700,136257/23200.6)
				};
				\addlegendentry{$Atc^\text{CSR}$}
				\addplot[
				color=teal,
				mark=triangle,
				]
				coordinates {
					(220,182.62/86.44)
					(740,536.79/212.01)
					(2692,2164.64/632.94)
					(10244,8316.04/1435.52)
					(39940,33425.2/4721.9)
					(157700,136257/18798.5)
				};
				\addlegendentry{$Sp^\text{CSR}$}
				\addplot[
				color=purple,
				mark=pentagon,
				]
				coordinates {
					(220,182.62/75.75)
					(740,536.79/165.82)
					(2692,2164.64/328.6)
					(10244,8316.04/730.97)
					(39940,33425.2/2611.03)
					(157700,136257/17787.2)
				};
				\addlegendentry{$Col^\text{CSR}_\text{vec}$}
				\addplot[
				color=violet,
				mark=square,
				]
				coordinates {
					(220,182.62/58.07)
					(740,536.79/132.39)
					(2692,2164.64/328.56)
					(10244,8316.04/752.74)
					(39940,33425.2/2010.53)
					(157700,136257/8301.37)
				};
				\addlegendentry{$Sp_\text{vec}^\text{CSR}$}
				\addplot[
				color=cyan,
				mark=o,
				]
				coordinates {
					(220,182.62/56.89)
					(740,536.79/128.22)
					(2692,2164.64/318.18)
					(10244,8316.04/712.11)
					(39940,33425.2/1807.95)
					(157700,136257/7471.05)
				};
				\addlegendentry{$Sp_\text{vec}^\text{CRAC}$}
			\end{axis}
		\end{tikzpicture}
		\caption{h-benchmark on unstructured meshes}
	\end{subfigure}
	\caption{Improvement factors of assembly times for the h-benchmark for the case $d = 4$.}
	\label{fig:par_hd4}
\end{figure}

Having considered the $d = 1$ case, we now take a look at benchmarks with $d = 4$, which translates to an alignment of 4 DOF on each node. The h-benchmark is depicted in \cref{fig:par_hd4}. The higher levels of continuity translate to the high performance of the vectorized methods. Interestingly, the colouring scheme suffers the most from the memory wall in terms of performance, whereas the other methods simply stagnate. 
Unlike in previous benchmarks, the CRAC format is consistently faster than the CSR format, peaking at an improvement of $13\%$. 

\begin{figure}
	\begin{subfigure}{0.48\linewidth}
		\centering
		\begin{tikzpicture}
			\pgfplotsset{every tick label/.append style={font=\normalsize},
				every x tick scale label/.append style={yshift=0.3em}}
			\begin{axis}[
				/pgf/number format/1000 sep={},
				axis lines = left,
				xlabel={Global degrees of freedom},
				ylabel={$c = t_\text{avg}^\text{seq} / t_\text{avg}$},
				xmin=0, xmax=150e3,
				ymin=0, ymax= 24,
				xtick={0, 30e3, 60e3, 90e3, 120e3, 150e3},
				ytick={0,4,8,12,16, 20,24},
				legend style={at={(0.5,0.03)},anchor=south},
				ymajorgrids=true,
				grid style=dotted,
				]
				\addplot[
				color=blue,
				mark=diamond,
				]
				coordinates {
					(9604,7463.98/1398.64)
					(37636,43359/7025.41)
					(84100,144142/24493.4)
					(148996,389784/62320.5)
				};
				\addlegendentry{$Atc^\text{CSR}$}
				\addplot[
				color=teal,
				mark=triangle,
				]
				coordinates {
					(9604,7463.98/1316.07)
					(37636,43359/5801.24)
					(84100,144142/19849)
					(148996,389784/51473)
				};
				\addlegendentry{$Sp^\text{CSR}$}
				\addplot[
				color=purple,
				mark=pentagon,
				]
				coordinates {
					(9604,7463.98/1009.34)
					(37636,43359/3839.82)
					(84100,144142/13249.6)
					(148996,389784/30511.8)
				};
				\addlegendentry{$Col^\text{CSR}_\text{vec}$}
				\addplot[
				color=violet,
				mark=square,
				]
				coordinates {
					(9604,7463.98/702.6)
					(37636,43359/2529.24)
					(84100,144142/10037.2)
					(148996,389784/25452.8)
				};
				\addlegendentry{$Sp_\text{vec}^\text{CSR}$}
				\addplot[
				color=cyan,
				mark=o,
				]
				coordinates {
					(9604,7463.98/649.23)
					(37636,43359/2029.98)
					(84100,144142/8608.14)
					(148996,389784/21081.3)
				};
				\addlegendentry{$Sp_\text{vec}^\text{CRAC}$}
			\end{axis}
		\end{tikzpicture}
		\caption{p-benchmark on structured meshes}
	\end{subfigure}
	\begin{subfigure}{0.48\linewidth}
		\centering
		\begin{tikzpicture}
			\pgfplotsset{every tick label/.append style={font=\normalsize},
				every x tick scale label/.append style={yshift=0.3em}}
			\begin{axis}[
				/pgf/number format/1000 sep={},
				axis lines = left,
				xlabel={Global degrees of freedom},
				ylabel={$c = t_\text{avg}^\text{seq} / t_\text{avg}$},
				xmin=0, xmax=160e3,
				ymin=0, ymax= 24,
				xtick={0, 40e3, 80e3, 120e3, 160e3},
				ytick={0,4,8,12,16,20,24},
				legend style={at={(0.5,0.03)},anchor=south},
				ymajorgrids=true,
				grid style=dotted,
				]
				\addplot[
				color=blue,
				mark=diamond,
				]
				coordinates {
					(10244,8262.94/1440.24)
					(39940,46409/7701.52)
					(89092,154510/25753.5)
					(157700,411737/66266.3)
				};
				\addlegendentry{$Atc^\text{CSR}$}
				\addplot[
				color=teal,
				mark=triangle,
				]
				coordinates {
					(10244,8262.94/1416.59)
					(39940,46409/6297.45)
					(89092,154510/21011.2)
					(157700,411737/61465.8)
				};
				\addlegendentry{$Sp^\text{CSR}$}
				\addplot[
				color=purple,
				mark=pentagon,
				]
				coordinates {
					(10244,8262.94/721.38)
					(39940,46409/4278.89)
					(89092,154510/13839.3)
					(157700,411737/32376.4)
				};
				\addlegendentry{$Col^\text{CSR}_\text{vec}$}
				\addplot[
				color=violet,
				mark=square,
				]
				coordinates {
					(10244,8262.94/753.76)
					(39940,46409/2748.4)
					(89092,154510/10312.7)
					(157700,411737/29008.3)
				};
				\addlegendentry{$Sp_\text{vec}^\text{CSR}$}
				\addplot[
				color=cyan,
				mark=o,
				]
				coordinates {
					(10244,8262.94/706.58)
					(39940,46409/2244.44)
					(89092,154510/8719.01)
					(157700,411737/24945.5)
				};
				\addlegendentry{$Sp_\text{vec}^\text{CRAC}$}
			\end{axis}
		\end{tikzpicture}
		\caption{p-benchmark on unstructured meshes}
	\end{subfigure}
	\caption{Improvement factors of assembly times for the p-benchmark for the case $d = 4$ using the fourth structured and unstructured meshes.}
	\label{fig:par_pd4}
\end{figure}
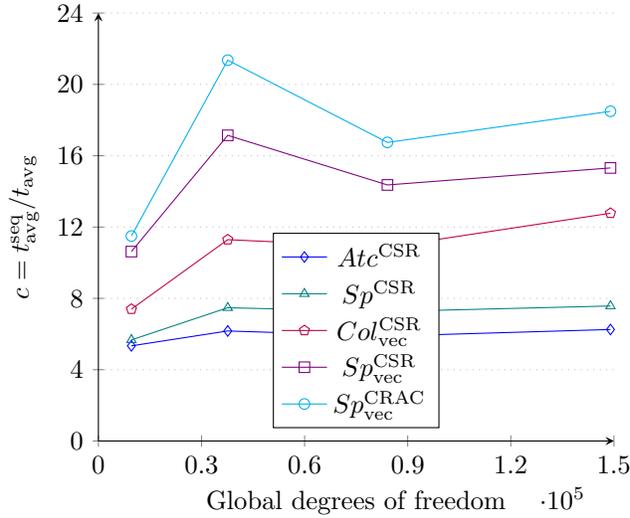
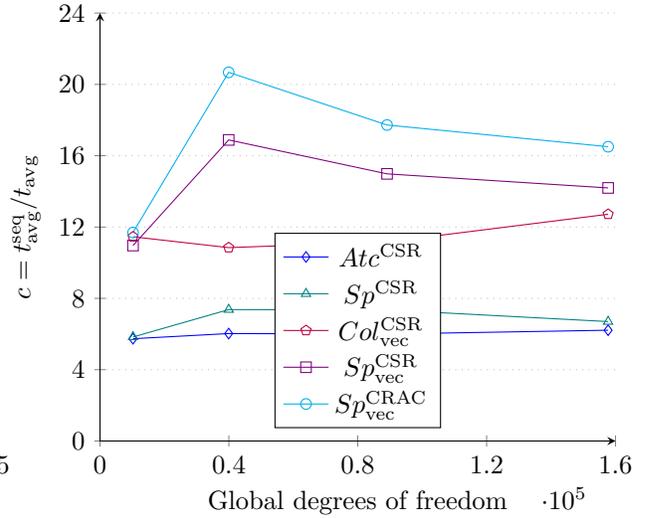

Next we consider the p-benchmark, see \cref{fig:par_pd4}, where the jump resulting from hitting the memory wall is apparent. Nevertheless, the performance improvement of the vectorized schemes is clear and once again, the CRAC format outperforms the CSR format.

\begin{figure}
	\begin{subfigure}{0.48\linewidth}
		\centering
		\begin{tikzpicture}
			\pgfplotsset{every tick label/.append style={font=\normalsize},
				every x tick scale label/.append style={yshift=0.3em}}
			\begin{axis}[
				/pgf/number format/1000 sep={},
				axis lines = left,
				xlabel={Global degrees of freedom},
				ylabel={$c = t_\text{avg}^\text{seq} / t_\text{avg}$},
				xmin=0, xmax=300e3,
				ymin=0, ymax= 28,
				xtick={0,50e3,100e3,150e3 ,200e3, 250e3, 300e3},
				ytick={0,4,8,12,16,20, 24, 28},
				legend pos=north west,
				ymajorgrids=true,
				grid style=dotted,
				]
				\addplot[
				color=blue,
				mark=diamond,
				]
				coordinates {
					(37249,7076.59/1361.01)
					(74498,31396.9/5034.61)
					(111747,69706.5/12502.5)
					(148996,128901/22629.1)
					(186245,209184/36242.6)
					(223494,305040/52876.7)
					(260743,372297/72190.4)
					(297992,557512/93593.8)
				};
				\addlegendentry{$Atc^\text{CSR}$}
				\addplot[
				color=teal,
				mark=triangle,
				]
				coordinates {
					(37249,7076.59/1442.46)
					(74498,31396.9/3898.17)
					(111747,69706.5/10199.3)
					(148996,128901/18332.4)
					(186245,209184/29245)
					(223494,305040/42127.4)
					(260743,372297/54387.7)
					(297992,557512/77242.6)
				};
				\addlegendentry{$Sp^\text{CSR}$}
				\addplot[
				color=purple,
				mark=pentagon,
				]
				coordinates {
					(37249,7076.59/1772.72)
					(74498,31396.9/2730.62)
					(111747,69706.5/8956.69)
					(148996,128901/16580.5)
					(186245,209184/27065.3)
					(223494,305040/38148.9)
					(260743,372297/51989.9)
					(297992,557512/58506.1)
				};
				\addlegendentry{$Col^\text{CSR}_\text{vec}$}
				\addplot[
				color=violet,
				mark=square,
				]
				coordinates {
					(37249,7076.59/1540.06)
					(74498,31396.9/2852.18)
					(111747,69706.5/5527.99)
					(148996,128901/8015.35)
					(186245,209184/11691.7)
					(223494,305040/15640.8)
					(260743,372297/20905.8)
					(297992,557512/26254.2)
				};
				\addlegendentry{$Sp_\text{vec}^\text{CSR}$}
				\addplot[
				color=cyan,
				mark=o,
				]
				coordinates {
					(37249,7076.59/1749.13)
					(74498,31396.9/2985.52)
					(111747,69706.5/5063.25)
					(148996,128901/7120.83)
					(186245,209184/9963.22)
					(223494,305040/13287.3)
					(260743,372297/17749.3)
					(297992,557512/21811.2)
				};
				\addlegendentry{$Sp_\text{vec}^\text{CRAC}$}
			\end{axis}
		\end{tikzpicture}
	\caption{d-benchmark on structured meshes}
	\end{subfigure}
    \begin{subfigure}{0.48\linewidth}
    	\centering
    	\begin{tikzpicture}
    		\pgfplotsset{every tick label/.append style={font=\normalsize},
    			every x tick scale label/.append style={yshift=0.3em}}
    		\begin{axis}[
    			/pgf/number format/1000 sep={},
    			axis lines = left,
    			xlabel={Global degrees of freedom},
    			ylabel={$c = t_\text{avg}^\text{seq} / t_\text{avg}$},
    			xmin=0, xmax=320e3,
    			ymin=0, ymax= 28,
    			xtick={0,80e3 ,160e3, 240e3, 320e3},
    			ytick={0,4,8,12,16,20, 24, 28},
    			legend pos=north west,
    			ymajorgrids=true,
    			grid style=dotted,
    			]
    			\addplot[
    			color=blue,
    			mark=diamond,
    			]
    			coordinates {
    				(39425,7523.68/1530.13)
    				(78850,34071.8/5467.05)
    				(118275,74094/13232.4)
    				(157700,136654/23617.8)
    				(197125,220003/37781.1)
    				(236550,319044/54280.1)
    				(275975,392885/72165.1)
    				(315400,580393/98640.1)
    			};
    			\addlegendentry{$Atc^\text{CSR}$}
    			\addplot[
    			color=teal,
    			mark=triangle,
    			]
    			coordinates {
    				(39425,7523.68/1566.46)
    				(78850,34071.8/4230.2)
    				(118275,74094/10748.9)
    				(157700,136654/19235.9)
    				(197125,220003/30448.5)
    				(236550,319044/43475.7)
    				(275975,392885/56910.9)
    				(315400,580393/80017.3)
    			};
    			\addlegendentry{$Sp^\text{CSR}$}
    			\addplot[
    			color=purple,
    			mark=pentagon,
    			]
    			coordinates {
    				(39425,7523.68/2156.48)
    				(78850,34071.8/3091.23)
    				(118275,74094/9733.77)
    				(157700,136654/17857.5)
    				(197125,220003/28827.8)
    				(236550,319044/40560.7)
    				(275975,392885/55148.4)
    				(315400,580393/61556.9)
    			};
    			\addlegendentry{$Col^\text{CSR}_\text{vec}$}
    			\addplot[
    			color=violet,
    			mark=square,
    			]
    			coordinates {
    				(39425,7523.68/1619.28)
    				(78850,34071.8/3081.8)
    				(118275,74094/5883.15)
    				(157700,136654/8413.59)
    				(197125,220003/12174.5)
    				(236550,319044/15931.2)
    				(275975,392885/21806.1)
    				(315400,580393/28017.4)
    			};
    			\addlegendentry{$Sp_\text{vec}^\text{CSR}$}
    			\addplot[
    			color=cyan,
    			mark=o,
    			]
    			coordinates {
    				(39425,7523.68/1925.77)
    				(78850,34071.8/3315.08)
    				(118275,74094/5438.96)
    				(157700,136654/7584.73)
    				(197125,220003/10469.3)
    				(236550,319044/13503.8)
    				(275975,392885/18457.4)
    				(315400,580393/23466.1)
    			};
    			\addlegendentry{$Sp_\text{vec}^\text{CRAC}$}
    		\end{axis}
    	\end{tikzpicture}
    	\caption{d-benchmark on unstructured meshes}
    \end{subfigure}
	\caption{Improvement factors of assembly times for the d-benchmark using the sixth structured and unstructured meshes.}
	\label{fig:par_d}
\end{figure}
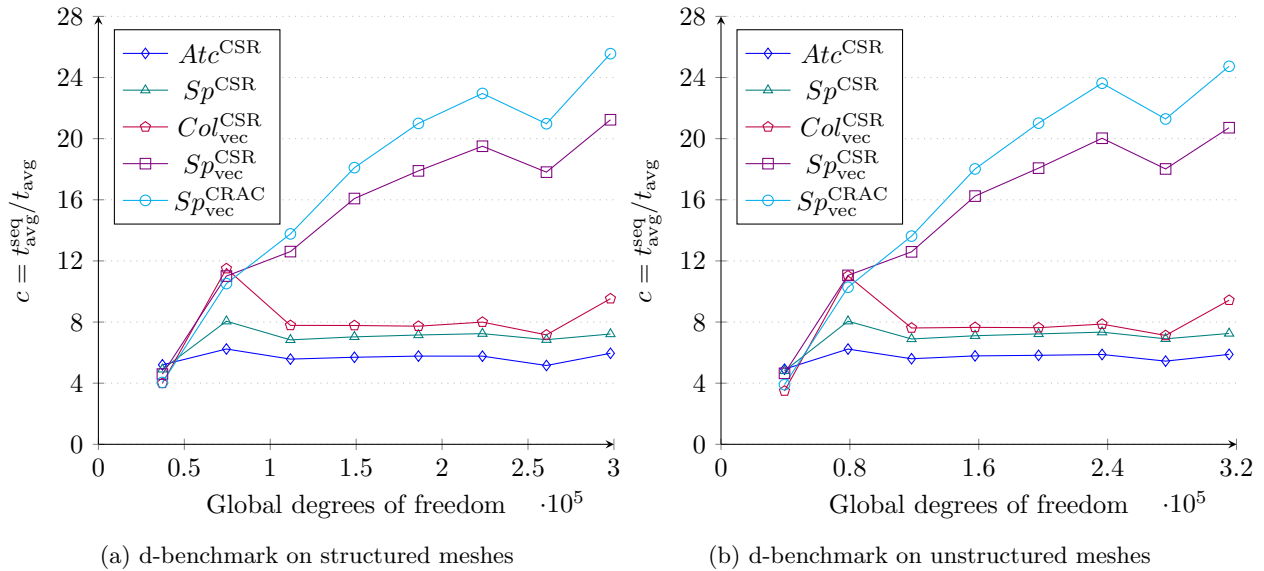

In the previous benchmark the dimension of the problem was set to $d=4$, which guaranteed a minimum alignment of 4 for all column sequences. The d-benchmark tests the performance improvement for an ever increasing alignment using the sixth meshes, see \cref{fig:par_d}. The L3 cache is surpassed for $d>2$. Even so, the vectorized \lstinline|spin_int| assembly schemes continue to improve, whereas the coloured scheme stagnates far earlier. Since also the coloured assembly scheme is vectorized, this suggests the method requires more frequent access to RAM. I.e., the coloured assembly scheme results in more cache-misses. The peak speed improvement factor of the \lstinline|spin_int| method with the CRAC format is $c = 25.56$ for the structured mesh, being $20\%$ faster than the CSR version. A comparison between the coloured scheme and the \lstinline|spin_int| version of the CSR format shows the latter to be $2.19-2.22$ times faster than the coloured scheme.
Another test with the fourth structured mesh and $d=20$ yields the improvement factor $c = 28.54$, which represents an approximate limit for highly aligned problems larger than the L3 cache for this assembly scheme on our system.

\subsection{Storage factor}
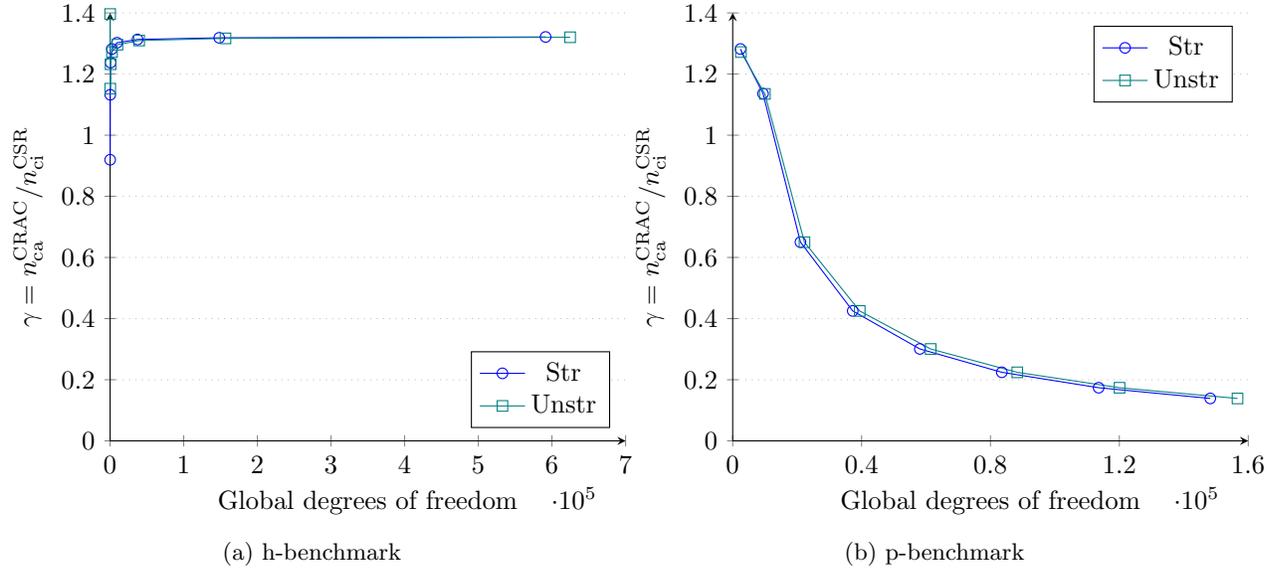
\begin{figure}
	\begin{subfigure}{0.48\linewidth}
		\centering
		\begin{tikzpicture}
			\pgfplotsset{every tick label/.append style={font=\normalsize},
				every x tick scale label/.append style={yshift=0.3em}}
			\begin{axis}[
				/pgf/number format/1000 sep={},
				axis lines = left,
				xlabel={Global degrees of freedom},
				ylabel={$\gamma = n^\text{CRAC}_\text{ca} / n^\text{CSR}_\text{ci} $},
				xmin=0, xmax=7e5,
				ymin=0, ymax= 1.4,
				xtick={0, 1e5, 2e5, 3e5, 4e5, 5e5, 6e5, 7e5},
				ytick={0,0.2,0.4,0.6,0.8,1,1.2,1.4},
				legend pos=south east,
				ymajorgrids=true,
				grid style=dotted,
				]
				\addplot[
				color=blue,
				mark=o,
				]
				coordinates {
					(49,332/361)
					(169,1550/1369)
					(625,6594/5329)
					(2401,26952/21025)
					(9409,108792/83521)
					(37249,437208/332929)
					(148225,1752984/1329409)
					(591361,7020312/5313025)
				};
				\addlegendentry{Str}
				\addplot[
				color=teal,
				mark=square,
				]
				coordinates {
					(55,546/391)
					(185,1688/1465)
					(673,6980/5665)
					(2561,28332/22273)
					(9985,114484/88321)
					(39425,460676/351745)
					(156673,1848612/1403905)
					(624641,7406692/5609473)
				};
				\addlegendentry{Unstr}
			\end{axis}
		\end{tikzpicture}
		\caption{h-benchmark}

	\end{subfigure}
	\begin{subfigure}{0.48\linewidth}
		\centering
		\begin{tikzpicture}
			\pgfplotsset{every tick label/.append style={font=\normalsize},
				every x tick scale label/.append style={yshift=0.3em}}
			\begin{axis}[
				/pgf/number format/1000 sep={},
				axis lines = left,
				xlabel={Global degrees of freedom},
				ylabel={$\gamma = n^\text{CRAC}_\text{ca} / n^\text{CSR}_\text{ci} $},
				xmin=0, xmax=160000,
				ymin=0, ymax=1.4,
				xtick={0,40000,80000,120000, 160000},
				ytick={0, 0.2, 0.4, 0.6, 0.8, 1, 1.2, 1.4},
				legend pos=north east,
				ymajorgrids=true,
				grid style=dotted,
				]
				\addplot[
				color=blue,
				mark=o,
				]
				coordinates {
					(2401,26952/21025)
					(9409,168378/148225)
					(21025,338052/519841)
					(37249,565198/1329409)
					(58081,849816/2825761)
					(83521,1191906/5313025)
					(113569,1591468/9150625)
					(148225,2048502/14753281)
				};
				\addlegendentry{Str}
				\addplot[
				color=teal,
				mark=square,
				]
				coordinates {
					(2561,28332/22273)
					(9985,177834/156673)
					(22273,356876/549121)
					(39425,596534/1403905)
					(61441,896808/2983681)
					(88321,1257698/5609473)
					(120065,1679204/9660673)
					(156673,2161326/15575041)
				};
				\addlegendentry{Unstr}
			\end{axis}
		\end{tikzpicture}
		\caption{p-benchmark}
	\end{subfigure}
	\caption{Column storage improvement factors for scalar problems $d = 1$.}
	\label{fig:hgamma}
\end{figure}

The speed improvement factors showed the CRAC format to be faster than the CSR format for problems with apparent alignment, whether the alignment occurs due to the dimension of the problem or the polynomial power of the element. Another factor for the efficiency of the CRAC of format is given by comparing the length of its column alignments array and the column indices array of the CSR format. \cref{fig:hgamma} depicts the behaviour of $\gamma$ for scalar-valued problems. For $d = 1$, h-adaption does not improve the factor, which is over $1$, suggesting the CSR format to be the lighter storage scheme. However, adapting $p$ to $3$ or higher results in higher alignment and the CRAC format becomes far lighter, having the lowest value of about $\gamma = 0.13$ for $p = 8$.
The improvement in storage for $p > 2$ can be explained by the way Gmsh \cite{eigenweb} defines nodes. Namely, first vertex nodes are numbered and edge nodes follow. For $p = 1$ no edge nodes are present. For $p = 2$ each edge contains one node, but they are not numbered contiguously with respect to the vertex nodes. In the case of $p = 3$ each edge incorporates two nodes, that are numbered contiguously, resulting in natural DOF alignment. 

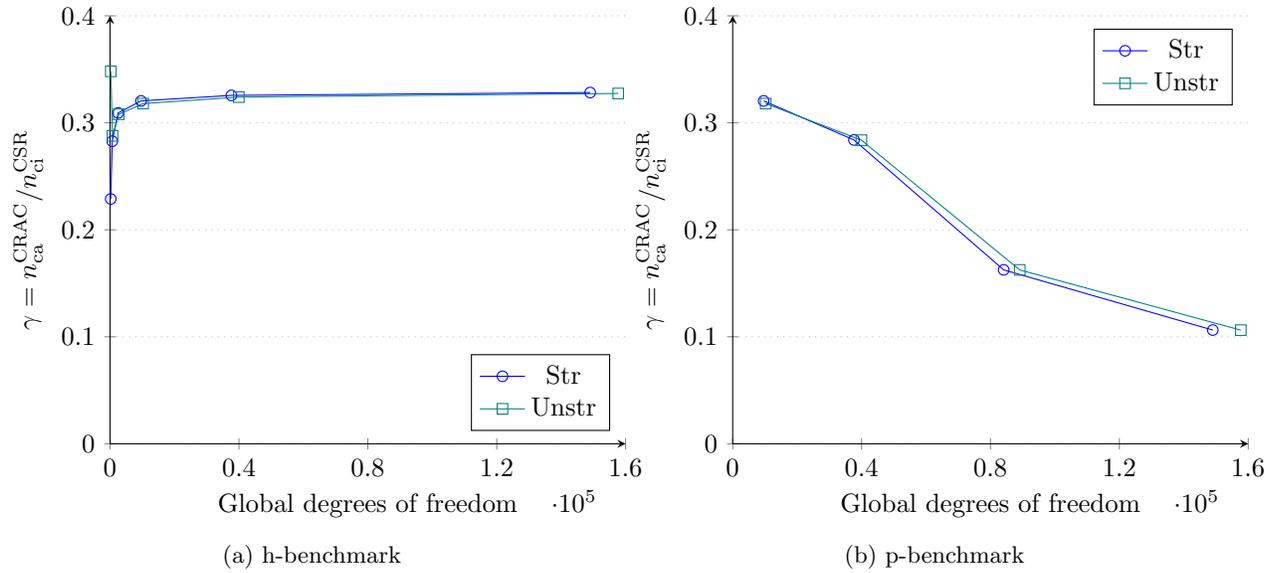
\begin{figure}
	\begin{subfigure}{0.48\linewidth}
		\centering
		\begin{tikzpicture}
			\pgfplotsset{every tick label/.append style={font=\normalsize},
				every x tick scale label/.append style={yshift=0.3em}}
			\begin{axis}[
				/pgf/number format/1000 sep={},
				axis lines = left,
				xlabel={Global degrees of freedom},
				ylabel={$\gamma = n^\text{CRAC}_\text{ca} / n^\text{CSR}_\text{ci} $},
				xmin=0, xmax=1.6e5,
				ymin=0, ymax= 0.4,
				xtick={0,0.4e5, 0.8e5, 1.2e5, 1.6e5},
				ytick={0, 0.1, 0.2, 0.3, 0.4},
				legend pos=south east,
				ymajorgrids=true,
				grid style=dotted,
				]
				\addplot[
				color=blue,
				mark=o,
				]
				coordinates {
					(196,1322/5776)
					(676,6194/21904)
					(2500,26370/85264)
					(9604,107802/336400)
					(37636,435162/1336336)
					(148996,1748826/5326864)
				};
				\addlegendentry{Str}
				\addplot[
				color=teal,
				mark=square,
				]
				coordinates {
					(220,2178/6256)
					(740,6746/23440)
					(2692,27914/90640)
					(10244,113322/356368)
					(39940,457930/1413136)
					(157700,1842698/5627920)
				};
				\addlegendentry{Unstr}
			\end{axis}
		\end{tikzpicture}
		\caption{h-benchmark}
	\end{subfigure}
	\begin{subfigure}{0.48\linewidth}
		\centering
		\begin{tikzpicture}
			\pgfplotsset{every tick label/.append style={font=\normalsize},
				every x tick scale label/.append style={yshift=0.3em}}
			\begin{axis}[
				/pgf/number format/1000 sep={},
				axis lines = left,
				xlabel={Global degrees of freedom},
				ylabel={$\gamma = n^\text{CRAC}_\text{ca} / n^\text{CSR}_\text{ci} $},
				xmin=0, xmax=160000,
				ymin=0, ymax=0.4,
				xtick={0,40000,80000,120000, 160000},
				ytick={0, 0.1, 0.2, 0.3, 0.4},
				legend pos=north east,
				ymajorgrids=true,
				grid style=dotted,
				]
				\addplot[
				color=blue,
				mark=o,
				]
				coordinates {
					(9604,107802/336400)
					(37636,673506/2371600)
					(84100,1352202/8317456)
					(148996,2260786/21270544)
				};
				\addlegendentry{Str}
				\addplot[
				color=teal,
				mark=square,
				]
				coordinates {
					(10244,113322/356368)
					(39940,711330/2506768)
					(89092,1427498/8785936)
					(157700,2386130/22462480)
				};
				\addlegendentry{Unstr}
			\end{axis}
		\end{tikzpicture}
		\caption{p-benchmark}
	\end{subfigure}
	\caption{Column storage improvement factors for scalar problems $d = 4$.}
	\label{fig:gammad4}
\end{figure}

The storage factor for the case $d = 4$ is given in \cref{fig:gammad4}. In the h-benchmark the factor remains approximately $\gamma = 0.32$. Using p-refinement allows to further lower the factor to $\gamma = 0.1$ for the same dimension of the problem.

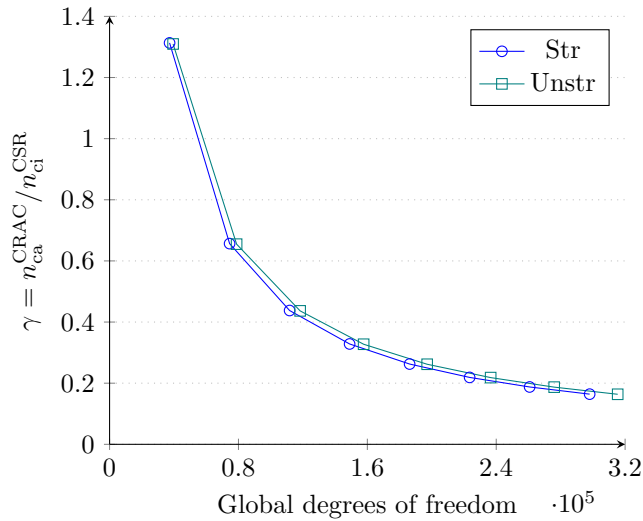
\begin{figure}
	\centering
	\begin{subfigure}{0.48\linewidth}
		\centering
		\begin{tikzpicture}
			\pgfplotsset{every tick label/.append style={font=\normalsize},
				every x tick scale label/.append style={yshift=0.3em}}
			\begin{axis}[
				/pgf/number format/1000 sep={},
				axis lines = left,
				xlabel={Global degrees of freedom},
				ylabel={$\gamma = n^\text{CRAC}_\text{ca} / n^\text{CSR}_\text{ci} $},
				xmin=0, xmax=320e3,
				ymin=0, ymax=1.4,
				xtick={0,80e3,160e3,240e3, 320e3},
				ytick={0, 0.2, 0.4, 0.6, 0.8, 1, 1.2, 1.4},
				legend pos=north east,
				ymajorgrids=true,
				grid style=dotted,
				]
				\addplot[
				color=blue,
				mark=o,
				]
				coordinates {
					(37249,437208/332929)
					(74498,874414/1331716)
					(111747,1311620/2996361)
					(148996,1748826/5326864)
					(186245,2186032/8323225)
					(223494,2623238/11985444)
					(260743,3060444/16313521)
					(297992,3497650/21307456)
				};
				\addlegendentry{Str}
				\addplot[
				color=teal,
				mark=square,
				]
				coordinates {
					(39425,460676/351745)
					(78850,921350/1406980)
					(118275,1382024/3165705)
					(157700,1842698/5627920)
					(197125,2303372/8793625)
					(236550,2764046/12662820)
					(275975,3224720/17235505)
					(315400,3685394/22511680)
				};
				\addlegendentry{Unstr}
			\end{axis}
		\end{tikzpicture}
	\end{subfigure}
	\caption{Column storage improvement factors for the d-benchmark.}
	\label{fig:gammad}
\end{figure}

Lastly, for higher dimensional problems the storage factor is given by \cref{fig:gammad}. Clearly, the CRAC format is the lighter storage scheme for $d > 1$.

\begin{remark}
		A lower $\gamma$ value translates to a lighter storage scheme and higher contiguity levels of the sparse matrix. The former results in better locality of the reference and the latter serves to increase vectorization levels.
\end{remark}

\section{Conclusions and outlook}\label{sec:out}
Although the colouring scheme allows for parallel lock-free assembly of the global sparse matrix, it does not outperform the \lstinline|spin_int| methods in our tests. This suggests that the colouring scheme leads to more cache misses and increased RAM access. The atomic method allows to modify every entry of the global sparse matrix in parallel but is also outperformed by the \lstinline|spin_int| methods. 
Consequently, it seems that the lock-unlock patterns of \lstinline|spin_int| methods do not dominate the assembly run time.

Using a format containing alignment information for the element \dofs array, we have been able to vectorize the assembly, yielding larger speed-up factors. This format is general and can be used to assemble all sorts of sparse matrices common in FE-Analysis. The introduced \lstinline|spin_int| data structures allow for data-race free assembly without increasing the storage requirements of the format and yield the best performance in their vectorized versions.

Based on the grouping of aligned columns, the CRAC format allows to improve the search algorithms needed for accessing coefficients, resulting in shorter assembly times.
Further, for the majority of the problems discussed, the CRAC format was both the lighter and faster storage scheme. Even for $d=1$ (only one DOF per node) the CRAC format leads to faster performance in case of p-adaption.

The CRAC format also entails other benefits not yet discussed in this paper. Since the format retains the alignment data of column indices for each row in the sparse matrix, this data can now be easily exploited for the design of parallel, vectorized versions of the Matrix-Vector product and correspondingly, the design of iterative solvers.
This information is not hard-coded into the sparse matrix design, but instead arises naturally from the variational problem and polynomial power of the elements, making the CRAC format a flexible and efficient design for dynamically defined FE-problems and the application of p- and hp-refinements. 

\bibliographystyle{spmpsci}      

\bibliography{ref}

\end{document}